\documentclass[12pt]{article}
\usepackage[english]{babel}
\usepackage{amsmath}
\usepackage{amssymb}

\newtheorem{theorem}{Theorem}[section]
\newtheorem{lemma}[theorem]{Lemma}

\newtheorem{remark}[theorem]{Remark}
\numberwithin{equation}{section}
\newcommand{\finedim}{\hfill\hbox{\enspace$\square$}\smallskip\\}
\newcommand{\iniziodim}{PROOF. }

\title{Sequential order under CH}

\author{Chiara Baldovino,\\ Department of Mathematics, University of Trieste\\
Via Alfonso Valerio 12/1, 34127 Trieste -- Italy,\\ chiara.baldovino@polito.it}

\begin{document}

\maketitle

\section{INTRODUCTION}
Let $X$ be a topological space and $M\subseteq X$;  the
\emph{sequential closure} of $M$ is $seqcl(M)=\{x\in X:\exists
(x_{n})_{n\in\omega}\subseteq M\footnote{By the notation
$(x_n)_{n\in\omega}\subset M$ we mean that $(x_n)_{n\in\omega}$ is
a sequence and that $x_n\in M$ for every $n\in\omega$.},
\lim_{n\in\omega}x_n=x \}.$ For every ordinal
$\alpha\leq\omega_1$, the \emph{$\alpha$-sequential closure} of
$M$ is inductively defined as follows:
\begin{itemize}
\item[-] $seqcl_0(M)=M$ and $seqcl_1(M)=seqcl(M)$;
\item[-] $seqcl_{\alpha+1}(M)=seqcl(seqcl_{\alpha}(M))$;
\item[-] $seqcl_{\alpha}(M)=\bigcup_{\beta<\alpha}seqcl_{\beta}(M)$ if $\alpha$ is a
limit ordinal.
\end{itemize}

A topological space $X$ is said to be \emph{sequential} if
$$seqcl_{\omega_1}(M)=\overline{M}, \qquad \forall M\subseteq X;$$ this definition is
equivalent to that according to which a space $X$ is said to be \emph{sequential} if
every sequentially closed subset of $X$ is closed.

The \emph{sequential order} of a sequential space $X$ is an ordinal invariant of the
space defined as $$\sigma(X)=\min\{\alpha\leq\omega_1:\forall M\subseteq X,
seqcl_{\alpha}(M)=\overline{M}\}.$$
 While the problem naturally posed in the sixties
concerning the possibility to produce examples of sequential spaces of any
sequential order up to and including $\omega_1$ in ZFC was completely solved in
the affirmative by Arhangel'ski\u{\i} and Franklin (cf. \cite{arhangelskii}),
it turns out difficult to construct compact sequential spaces without
additional assumptions of the Theory of the Sets, even of sequential order 3;
indeed, up to now, $2$ is the maximum order of sequentiality of a compact space
in ZFC. In this context the work due to Ba\v{s}kirov and concisely presented in
a Doklady article (see \cite{baskirov}) gathers a certain prominence: in this
paper the author suggests a scheme of construction to produce compact
sequential spaces of any order as quotient spaces of $\beta\omega$ under the
assumption of the Continuum Hypothesis. Since Ba\v{s}kirov's work is very
concise and devoid of any proof and check, we write down the construction with
some essential alterations with respect to the original work in order to
complete it in all details and explain
where the Continuum Hypothesis is essentially used.\\
In order to make comprehensive our survey concerning the
possibility to obtain compact sequential spaces of the greatest
order, we have to mention the construction under CH due to Kannan
(cf.
 \cite{kannan}) and the more recent constructions under MA
  due to Dow (see \cite{dow1} and \cite{dow2}).
Under the Continuum Hypothesis, Kannan ensures that it is possible
to construct compact sequential spaces of any order while, under
MA, Dow manages to give an example of a compact sequential space
of
order $4$, the best upper bounds under this axiom up to now.\\
 While Ba\v{s}kirov
suggests a construction from top to down, Kannan and Dow present a construction
from down to top. Indeed Ba\v{s}kirov works in $\beta\omega$ and by assuming to
have constructed all the spaces of sequential order a successor ordinal less
than a fixed successor ordinal $\alpha+1$ he gives a starting decomposition on
$\beta\omega$; the Continuum Hypothesis guarantees him that in $\omega_1$ steps
he can purify the starting decomposition in such a way that in the space
associated to the last decomposition there is a new point fit to produce a
space of sequential order $\alpha+1$. Instead Kannan and Dow start from the
natural numbers with the discrete topology. If we want to summarize the idea of
Kannan, we can say that he generalizes the construction of the one-point
compactification of the Mr\'{o}wka-Isbell space. On the other hand, Dow
constructs by transfinite induction on $\mathfrak{c}$ three suitable families
of subsets of $\omega$ in such a way that the Stone space associated to the
Boolean algebra generated by the elements of these subsets admits a point of
sequential order $4$.\\
There is a remarkable reason to determine the maximum possible sequential order
in the presence of the PFA which implies Martin's axiom and
$\mathfrak{c}=\omega_2$; indeed in 1989 Balogh solved the Moore-Mr\'{o}wka
problem proving that each compact space of countable tightness is sequential
under PFA (see \cite{balogh}). If there is some finite bound on the sequential
order of compact sequential spaces in models of PFA, it would mean that compact
spaces of countable tightness are a few steps away from being
Fr\'{e}chet-Urysohn. In \cite[Proposition~3.1]{dow2}, Dow points out that there
are obstructions to extend his type of construction to produce compact
sequential spaces of order greater than $4$. However the problem if there
exists a bound on the sequential order of compact sequential spaces in models
of PFA is still open.

\section{PRELIMINARY FACTS}

 We are interested in the construction of compact sequential spaces
of sequential order $1$ and $2$ in ZFC as quotient spaces of
$\beta\omega$; we choose to present the following constructions to
enter into the scheme of the
main construction we will present.\\
 Let us consider the
space
$$K_{1}=\beta\omega/\omega^{\ast}=\beta\omega/\approx_{1}$$ where
$$x\approx_{1} y\Leftrightarrow (x=y \vee (x\in{\omega^{\ast}}\wedge
y\in{\omega^{\ast}}))$$ and let us denote by $j_{1}$ the natural
quotient mapping from $\beta\omega$ to $K_{1}$.
 Trivially the
one-point elements of the quotient under the relation of equivalence
$\approx_{1}$ are images of the points of $\omega$, while the natural quotient
mapping collapses all the free ultrafilters to a single point $P$.
 The topology of the space
$K_{1}=\beta\omega/\omega^{\ast}=\omega\cup\left\{P\right\}$ with
$P\notin{\omega}$ is a topology $\tau$ such that the points of
$\omega$ are isolated while $P$ has a fundamental system of (open)
neighborhoods formed by $\left\{U_{P,F}
\right\}=\left\{\left\{P\right\}\cup j_{1}(\omega\backslash F):
F\in[\omega]^{<\omega}\right\}$ as we prove in the following
lemma.

\begin{lemma}\label{lemma intorni punti in K1} Let $P=j_{1}(\omega^{\ast})\in K_{1}$;
a fundamental system of neighborhoods of $P$ is given by the collection
$\left\{j_{1}(\beta\omega\backslash F): F\in[\omega]^{<\omega}\right\}$.
\end{lemma}
\iniziodim If $A=\beta\omega\backslash F$ with
$F\in{\left[\omega\right]}^{<\omega}$,
 then $A$ is a saturated open subset of $\beta\omega$ and $A\supseteq\omega^{\ast}$
  whence $j_{1}(A)$ is an open neighborhood of
 $P$.\\
Suppose now that $V$ is an arbitrary neighborhood of $P$ in $K_{1}$: we want to
prove that there exists $F\in [\omega]^{<\omega}$ such that
$j_{1}(\beta\omega\backslash F)\subseteq V$, i.e. such that
$\beta\omega\backslash F\subseteq j_{1}^{-1}(V)$. Towards a contradiction
suppose that $|\beta\omega\backslash j_{1}^{-1}(V)|=\omega$; then we can set
$L=\beta\omega\backslash j_{1}^{-1}(V)$ and consider a free ultrafilter
$\mathcal{U}\in\omega^{\ast}\subseteq j_{1}^{-1}(V)$ such that $L\in
\mathcal{U}$. Since $j_{1}^{-1}(V)$ is open, there exists an infinite set $H$
of $\omega$ with $H\in \mathcal{U}$ such that $H^{\ast}\cup H\subseteq
j_{1}^{-1}(V)$ (in particular $H\subseteq j_{1}^{-1}(V)$). Now $H,L\in
\mathcal{U}$ and then it follows that $H\cap L\neq \emptyset$ (and still better
$|H\cap L|=\omega$). Hence we can fix a point $m\in H\cap L$: on the one hand,
it holds that $m\in H\subseteq j_{1}^{-1}(V)$, while, on the other hand, it
results that $m\in L\subseteq \beta\omega\backslash j_{1}^{-1}(V)$. A
contradiction. \finedim

Notice that the fundamental neighborhoods $\{U_{P,F}\}$ of $P$ are clopen
subsets; we refer to these neighborhoods as \emph{elementary}.\\ It is easy to
see that the space $K_{1}$ has the same topology as a convergent sequence and
hence it is trivially a Hausdorff compact space.

\vspace{0,3cm}

Now we want to find a suitable relation of equivalence in $\beta\omega$ in such
a way that the corresponding quotient space is a Hausdorff compact space of
sequential order $2$. Let $\mathcal{M}$ be an infinite MAD family on $\omega$;
to every element $M\in\mathcal{M}$ we can associate the unique element
$M^{\ast}\subseteq \omega^{\ast}$ in the following way:
\begin{displaymath}
M  \mapsto M^{\ast}=\{\mathcal{U}\in{\omega^{\ast}}:M\in{\mathcal{U}}\}.
\end{displaymath}
It turns out that if $M_{1},M_{2}\in \mathcal{M}$ with $M_{1}\neq M_{2}$ then
$M_{1}^{\ast}\cap M_{2}^{\ast}=\emptyset$: indeed if the intersection was not
empty, then there would exist a free ultrafilter $\mathcal{V}\in\omega^{\ast}$
such that $M_{1}\in\mathcal{V}$ and $M_{2}\in\mathcal{V}$; thus $M_{1}\cap
M_{2}\in\mathcal{V}$ but $|M_{1}\cap M_{2}|<\omega$ and the ultrafilter would
be fixed against the hypothesis.\\
We also remark that the subset
$\omega^{\ast}\backslash\bigcup_{M\in{\mathcal{M}}}M^{\ast}$ is not empty: if
it was empty, then $\{M^{\ast}\cup \omega:M\in{\mathcal{M}}\}$ would be an
infinite and open cover of $\beta\omega$ from which it would be impossible to
extract a finite subcover; this fact clashes with the compactness of
$\beta\omega$.

Let us take into account the space $K_{2}=\beta\omega/\approx_{2}$ where two
free ultrafilter of $\beta\omega$ are equivalent under the relation
$\approx_{2}$ if they belong to the same $M^{\ast}$ with $M\in\mathcal{M}$;
moreover all the free ultrafilters belonging to the subset
$\omega^{\ast}\backslash\bigcup_{M\in{\mathcal{M}}}M^{\ast}$ are equivalent,
while the relation does not identify any point in $\omega$. Then, if we denote
by $j_{2}$ the natural quotient mapping from $\beta\omega$ to $K_{2}$, it holds
that $j_{2}$ leaves the points of $\omega$ unaltered, while it collapses every
$M^{\ast}$ with $M\in \mathcal{M}$ to a single point and the non-empty subset
$\omega^{\ast}\backslash\bigcup_{M\in{\mathcal{M}}}M^{\ast}$ to another single
point too.\\ Let us set $L_{0}=j_{2}(\omega)$,
$L_{1}=\{j_{2}(M^{\ast}):M\in\mathcal{M}\}$ and
$x_{\infty}=j_{2}(\omega^{\ast}\backslash\bigcup_{M\in{\mathcal{M}}}M^{\ast})$.
We say that the points in the set $L_{0}$ have level $0$ while the points of
$L_{1}$ have level $1$ and the point $x_{\infty}$ has level $2$. We will prove
that the levels of the points coincide with their sequential order with respect
to the set $L_{0}$.\\ We already know that the points in $L_{0}$ are isolated
in $K_{2}$; now we prove the following lemma about the neighborhoods of the
points of level $1$.
\begin{lemma}\label{lemma intorni punti ord1 in K2} Let $y\in L_{1}$ with $y=j_{2}(M^{\ast})$ and $M\in\mathcal{M}$;
 then the collection $\{U_{y,F}\}=\{\{y\}\cup
j_{2}(M\backslash F):F\in{[M]^{<\omega}}\}$ is a fundamental
system of neighborhoods for $y$.
\end{lemma}
\iniziodim On the one hand, since $M^{\ast}\cup{M}$ is a saturated open subset
of $\beta\omega$ and since every fixed ultrafilter in $\beta\omega$ is
isolated, it is evident that $M^{\ast}\cup(M\backslash F)$ is a saturated open
subset of $\beta\omega$ too for every $F\in{[M]^{<\omega}}$. Therefore
$j_{2}(M^{\ast}\cup(M\backslash F))=\{y\}\cup j_2(M\backslash F)$ is an open
neighborhood of $y$ for every $F\in[M]^{<\omega}$.
\\ On the other hand, let us suppose that $V$ is an arbitrary
neighborhood of $y$; we want to prove that there exists $F\in [M]^{<\omega}$
such that $j_{2}(M\backslash F)\subseteq V$, i.e. such that $M\backslash
F\subseteq j_{2}^{-1}(V)$. Let us suppose by contradiction that $|M\backslash
j_{2}^{-1}(V)|=\omega$; then let us set $L=M\backslash j_{2}^{-1}(V)$ and
consider a free ultrafilter $\mathcal{U}$ with $L\in \mathcal{U}$: it holds
that $M\in \mathcal{U}$ (since $L\subseteq M$) and then that $\mathcal{U}\in
M^{\ast}\subseteq j_{2}^{-1}(V)$. As $j_{2}^{-1}(V)$ is open, there exists
$H\in \mathcal{U}$ such that $H^{\ast}\cup H\subseteq j_{2}^{-1}(V)$ (in
particular $H\subseteq j_{2}^{-1}(V)$). In this way it turns out that $H,L\in
\mathcal{U}$; then $H\cap L\in\mathcal{U}$ whence $H\cap L\neq \emptyset$ (even
better $|H\cap L|=\omega$ since $\mathcal{U}$ is free). Let us fix a point
$m\in H\cap L$: on the one hand, it holds that $m\in H\subseteq j_{2}^{-1}(V)$,
while, on the other hand, it results that $m\in L\subseteq M\backslash
j_{2}^{-1}(V)$. A contradiction.
 \finedim

We want to remark that the fundamental neighborhoods $U_{y,F}$ of
a general point $y$ of level $1$ are clopen in $K_{2}$: indeed
$M^{\ast}\cup (M\backslash F)$ is a saturated closed subset of
$\beta\omega$ for each $F\in{[M]^{<\omega}}$; let us call
\emph{elementary} these neighborhoods.\\ We state in advance the
following remark to the lemma about
  the neighborhoods of the point $x_{\infty}$.

\begin{remark}\label{oss int ord2 in K2} For every $D\in[\omega]^{\omega}$
the subfamily $\mathcal{M}_{D}=\{M\in\mathcal{M}:|M\cap
D|=\omega\}$ is such that
$\bigcup\mathcal{M}_{D}\supseteq^{\ast}D$; suppose by
contradiction that $\bigcup\mathcal{M}_{D}\nsupseteq^{\ast}D$,
i.e. suppose that there exists a subset $E\subset D$ such that
$|E|=\omega$ and $(\bigcup\mathcal{M}_{D})\cap E=^{*}\emptyset$.
Since the family $\mathcal{M}$ is maximal, there exists a subset
$A\in\mathcal{M}\backslash\mathcal{M}_{D}$ such that $|A\cap
E|=\omega$ but $A\cap E\subset D$ and hence $|A\cap D|=\omega$;
therefore $A$ is an
element of $\mathcal{M}_{D}$ but this contradicts what we have just supposed.\\
Let us fix an arbitrary $D\in[\omega]^{\omega}$; if for every
finite union $\bigcup_{M\in \mathcal{F}}M$ with
$\mathcal{F}\in[\mathcal{M}]^{<\omega}$ it turns out that
$\bigcup_{M\in \mathcal{F}}M\nsupseteq^{\ast}D$, then
$|\mathcal{M}_{D}|\geq \omega$. Indeed if
$|\mathcal{M}_{D}|<\omega$, by taking
$\mathcal{F}=\mathcal{M}_{D}$, it holds that
$\bigcup\mathcal{F}\supseteq^{\ast}D$ because of what we have just
remarked above.
\end{remark}

  \begin{lemma}\label{lemma intorni punti ord2 in K2} The collection of the clopen subsets $K_{2}\backslash
   \bigcup_{x\in G}U_{x}$ (where $G$ is a finite set and for every $x\in G$ the clopen subset $U_{x}$
 is an elementary neighborhood of the point $x$ in $K_{2}$
 that can have level $0$ or $1$)
  is a base at the point $x_{\infty}$.

\end{lemma}
\iniziodim In an obvious way, $K_{2}\backslash \bigcup_{x\in G}U_{x}$ (where
$G$ is a finite set and for every $x\in G$ the clopen subset $U_{x}$
 is an elementary neighborhood of the point $x$ of level $0$ or $1$) is open and closed in $K_{2}$ since its complementary subset is a finite
  union of clopen subsets.\\
Now let $A$ be an open subset containing $x_{\infty}$ and let
$C=K_{2}\backslash A$ be the complementary closed subset. For every $x\in C$,
let $U_{x}$ be an elementary clopen neighborhood of $x$; trivially, by taking
all the clopen subsets $U_{x}$ with $x\in C$, we cover $C$. Let us consider
$j_{2}^{-1}(C)$: it is a closed subset in $\beta\omega$ and then it is compact.
The subsets $j_{2}^{-1}(U_{x})$ (with $x\in C$) form an open cover of
$j_{2}^{-1}(C)$; then there exists a finite subcover $\bigcup_{x\in
G}j_{2}^{-1}(U_{x})\supseteq j_{2}^{-1}(C)$. It turns out that
$$j_{2}(\bigcup_{x\in G}j_{2}^{-1}(U_{x}))=\bigcup_{x\in
G}j_{2}(j_{2}^{-1}(U_{x}))=\bigcup_{x\in G}U_{x}\supseteq
j_{2}(j_{2}^{-1}(C))=C;$$ by returning to the complementary subsets, we are
able to conclude that $K_{2}\backslash \bigcup_{x\in G}U_{x}\subseteq
K_{2}\backslash C=A.$ \finedim

 Let us call \emph{elementary} these clopen neighborhoods of
the point $x_{\infty}$. Now we are finally able to prove the
following lemma.
\begin{lemma}
$K_{2}$ is a compact sequential Hausdorff space of sequential
order $2$.
\end{lemma}
\iniziodim It is trivial to part any point $j_{2}(n)\in L_{0}$ from any other
point $Q\in K_{2}$: indeed we can take respectively the open disjoint
neighborhoods $j_{2}(\{n\})$ and $j_2(\beta\omega\backslash\{n\})$. We have to
analyse the other following two cases in order to conclude that $K_2$ is a
Hausdorff space
\begin{enumerate}
\item[1.] If we have to separate two points $P_{1},P_{2}\in K_{2}$
of level $1$, i.e. such that $j_2^{-1}(P_{1})=M_{1}^{\ast}$,
$j_2^{-1}(P_{2})=M_{2}^{\ast}$ with $M_{1},M_{2}\in{\mathcal{M}}$, then we
notice that $M_{1}^{\ast}\cup M_{1}$ and $M_{2}^{\ast}\cup M_{2}$ are open
subsets of $\beta\omega$ and that $F=M_{1}\cap M_{2}$ is finite and hence
closed in $\beta\omega$; thus it turns out that $M_{1}^{\ast}\cup
(M_{1}\backslash F)$
  and $M_{2}^{\ast}\cup (M_{2}\backslash F)$
  are disjoint saturated open subset of $\beta\omega$ and then their
   images $j_2(M_{1}^{\ast}\cup (M_{1}\backslash F))$
and $j_2(M_{2}^{\ast}\cup (M_{2}\backslash F))$ are disjoint open neighborhoods of
$P_{1}$ and $P_{2}$ respectively.
 \item[2.] If we have to part the point $x_{\infty}=j_{2}(\omega^{\ast}
 \backslash\bigcup_{M\in{\mathcal{M}}}M^{\ast})$
  from any point $P$
  of level $1$, i.e. such that $j_2^{-1}(P)=M_{1}^{\ast}$ with $M_1\in\mathcal{M}$, then we notice that
   $M_{1}^{\ast}\cup M_{1}$ is a saturated clopen subset of $\beta\omega$
    and that also its complement is saturated and clopen; therefore $j_2(M_{1}^{\ast}\cup M_{1})$
    and
$j_2(\beta\omega\backslash(M_{1}^{\ast}\cup M_{1}))$ are disjoint open neighborhoods of
$P$ and $x_{\infty}$ respectively (the latter one is a neighborhood of $x_{\infty}$ since
$$j_2^{-1}(x_{\infty})=\omega^{\ast}\backslash
\bigcup_{M\in{\mathcal{M}}}M^{\ast}\subseteq \beta\omega\backslash(M_{1}^{\ast}\cup
M_{1}).$$
\end{enumerate}
We can also conclude that $K_{2}$ is compact, since $j_{2}$ is a continuous
function from the compact space $\beta\omega$ to the Hausdorff space $K_{2}$.\\
Finally we can prove that the space $K_{2}$ has sequential order $2$; more
precisely we show that its sequential order, $\sigma(K_{2})$, is less than or
equal to $2$ (and, in particular,
that $K_2$ is a sequential space) and then that it is exactly $2$.\\
Let $S\subseteq K_{2}$ be an arbitrary subset such that $\overline{S}\neq{S}$:
we want to show that $\overline{S}=\emph{seqcl}_{2}(S)$; clearly it is enough
to prove that $\overline{S}\backslash S\subseteq\emph{seqcl}_{2}(S)$. Let us
consider a point $x\in{\overline{S}\backslash S}$ and let us set
$P=j_2^{-1}(x)$. Trivially $P\notin\{\{n\}:n\in{\omega}\}$: indeed the images
of points of $\omega$ in ${K_{2}}$ can not belong to $\overline{S}\backslash S$
since they are isolated. Therefore there are two cases to study:
\begin{enumerate}
\item [$1^{st}$)] $P=M^{\ast}$ with $M\in{\mathcal{M}}$; \item [$2^{nd}$)]
$P=\omega^{\ast}\backslash(\bigcup_{M\in{\mathcal{M}}}M^{\ast})$.
\end{enumerate}
Let us take them into account:
\begin{enumerate}
\item [$1^{st}$)] First of all if $x\in{\overline{S}\backslash S}$, by
Lemma \ref{lemma intorni punti ord1 in K2} we can assert that $x\in{\overline{S\cap
j_2(M)}}$. This fact allows us to conclude that $S\cap j_2(M)$ is infinite, since $K_{2}$
satisfies the $T_{1}$ separation axiom, and then we can write $S\cap j_2(M)$ as
$\{j_2(m_{n}): n\in{\omega}\}$ where $n\mapsto m_{n}$ is injective and $m_{n}\in{M}$ for
every $n\in{\omega}$. We prove that, for every neighborhood $V$ of $x$ in $K_{2}$,
$j_2(M)\backslash V$ is finite - this will imply that
$(j_2(m_{n}))_{n\in{\omega}}\rightarrow x$ in $K_{2}$ and hence that
$x\in{\emph{seqcl}_{1}(S)\subseteq{\emph{seqcl}_{2}(S)}}$. Towards a contradiction,
suppose that $j_2(M)\backslash V$ is infinite; then there exists $M'\in{[M]^{\omega}}$
such that $j_2(M')=j_2(M)\backslash V$. Let us consider a free ultrafilter $\mathcal{U}$
with $M'\in{\mathcal{U}}$. Since $V$ is an open neighborhood of $x$ in $K_{2}$, it turns
out that $j_2^{-1}(V)$ is a saturated open subset of $\beta\omega$ such that
$j_2^{-1}(x)=M^{\ast}\subseteq{j_2^{-1}(V)}$. Then it holds that
$\mathcal{U}\in{(M')^{\ast}}\subseteq{M^{\ast}}\subseteq{j_2^{-1}(V)}$ and hence there
exists an infinite subset $T\in{\mathcal{U}}$ such that $T^{\ast}\cup T\subseteq
j_2^{-1}(V)$ and, in particular, it turns out that $j_2(T)\subseteq j_2(T\cup
T^{\ast})\subseteq{V}$; since $T\cap{M'}\neq\emptyset$ (they are both elements of
$\mathcal{U}$), we can fix $h\in{T\cap{M'}}$ and on the one hand we obtain that
$j_2(h)\in{j_2(T)}\subseteq{V}$, while on the other hand we have
$j_2(h)\in{j_2(M')}=j_2(M)\backslash V$ and we reach a contradiction.

\item [$2^{nd}$)] Since $x\in\overline{S}$ by Lemma \ref{lemma
intorni punti ord2 in K2} it holds that that either in $S$ there are at least a
countable infinity of points $y_{n}\in L_{1}$ or, if $|S\cap L_{1}|<\infty$, in
$S$ there are infinite points of level $0$ (let us call $D$ the set consisting
of these points) such that it is not possible to cover $D$ with a finite number
of elementary neighborhoods of points of level $1$. In the former case by Lemma
\ref{lemma intorni punti ord2 in K2} any sequence extracted from $S\cap L_{1}$
converges to $x_{\infty}$ and then $x\in seqcl_{1}(S)$. In the latter case by
Remark \ref{lemma intorni punti ord2 in K2}, since it is not possible to cover
$D$ with a finite number of elementary neighborhoods of points of level $1$,
i.e. since for every finite union of elements $M\in\mathcal{M}$ it turns out
that $\bigcup M\nsupseteq^{\ast}D$, it follows that $|\mathcal{M}_{D}|=\infty$:
then there are infinite points of level $1$, $\{y_{\alpha}\}_{\alpha\in
A}\subseteq L_{1}$, such that $\{y_{\alpha}\}_{\alpha\in A}\subseteq
seqcl_{1}(j_{2}(D))\subseteq seqcl_{1}(S)$ by the former case; moreover a
sequence extracted from this set has to converge to $x_{\infty}$ as we have
remarked above and hence it turns out that $x_{\infty}\in seqcl_{2}(S)$.
\end{enumerate}
Now we have to prove that the sequential order of $K_{2}$ is exactly $2$: we
will show that there is a subset $D\subseteq K_{2}$ with the property that
$x_{\infty}\in\overline{D}$ and that no sequence extracted from $D$ converges
to $x_{\infty}$. Let us consider the subset $j_{2}(\omega)$; it is clear that
$x_{\infty}\in\overline{j_{2}(\omega)}$. Now we prove that no sequence
extracted from $j_2(\omega)$ converges to $x_{\infty}$; towards a
contradiction, suppose that there exists a sequence
$(j_2(m_{n}))_{n\in{\omega}}\rightarrow x_{\infty}$. In an obvious way, the set
$H=\{m_{n}:n\in{\omega}\}$ is infinite since $K_{2}$ satisfies the $T_{1}$
separation axiom; hence there exists an infinite set
$\tilde{M}\in{\mathcal{M}}$ such that $|H\cap\tilde{M}|=\omega$. If we call
$\tilde{y}$ the unique point in $L_{1}$ such that
$j_2^{-1}(\tilde{y})=\tilde{M}^{\ast}$, by Lemma \ref{lemma intorni punti ord1
in K2} it turns out that $\tilde{y}\in \overline{j_2(H)}$. Therefore we are in
the first situation we have studied above and so we can assert that there is a
subsequence $(m_{n_i})_{i\in{\omega}}$ of $(m_{n})_{n\in{\omega}}$ such that
$(j_2(m_{n_{i}}))_{i\in{\omega}}\rightarrow \tilde{y}$; on the other hand by
hypothesis we know that $(j_2(m_{n}))_{n\in{\omega}}\rightarrow x_{\infty}$ and
hence that $(j_2(m_{n_{i}}))_{i\in{\omega}}\rightarrow x_{\infty}$. This is a
contradiction since $K_{2}$ is a Hausdorff space and
$x_{\infty}\neq{\tilde{y}}$.  \finedim

We want to remark that the space $K_2$ we have just constructed is trivially
homeomorphic to the one-point compactification of the Mr\'{o}wka-Isbell space
$\Psi(\mathcal{M})$.\\ By referring to the type of construction of the space
$K_{2}$, we could think that a good idea to construct a space with a larger
order of sequentiality could be to associate a new infinite MAD family
$\mathcal{H}_{M}$ to every $M\in \mathcal{M}$; we will prove that in this way
we do not construct a space of higher sequential order.\\ Then let
$\mathcal{M}$ be an infinite MAD family on $\omega$ and let us suppose to
associate a new infinite MAD family $\mathcal{H}_{M}$ to every $M\in
\mathcal{M}$; let us consider the partition $\mathcal{P}=\{\{n\}:n\in
\omega\}\cup\{H^{\ast}:H\in \bigcup_{M\in
\mathcal{M}}\mathcal{H}_{M}\}\cup\{M^{\ast}\backslash \bigcup_{H\in
{\mathcal{H}_{M}}}H^{\ast}:M\in
{\mathcal{M}}\}\cup\{\omega^{\ast}\backslash\bigcup_{M\in\mathcal{M}}M^{\ast}\}$.
Let us set $K_{\textrm{\scriptsize{III}}}=\beta\omega/\approx$ where $\approx$
is the relation of equivalence associated to $\mathcal{P}$ and let
$j_{\textrm{\scriptsize{III}}}$ be the natural quotient mapping from
$\beta\omega$ to $K_{\textrm{\scriptsize{III}}}$. The space
$K_{\textrm{\scriptsize{III}}}$ consists of the following elements.
\begin{enumerate}
\item[$1^{st}$)] The isolated points of the form
$j_{\textrm{\scriptsize{III}}}(n)$ for every $n\in\omega$.
\item[$2^{nd}$)] The points of the form
$x_{H}=j_{\textrm{\scriptsize{III}}}(H^{\ast})$ for every $H\in\mathcal{H}_{M}$ and every
$M\in \mathcal{M}$; a fundamental system of neighborhoods for $x_H$ is given by
$$\{\{x_H\}\cup j_{\textrm{\scriptsize{III}}}(H\backslash
F)\}_{F\in{[H]^{<\omega}}}.$$ We refer to these neighborhoods with the symbols
$U_{x_{H},F}$. \item[$3^{rd}$)] The points of the form
$y_{M}=j_{\textrm{\scriptsize{III}}}(M^{\ast}\backslash \bigcup_{H\in
{\mathcal{H}_{M}}}H^{\ast})$ for every $M\in \mathcal{M}$; a fundamental system of
neighborhoods for $y_M$ is given by
$$\{\left(\{y_M\}\cup j_{\textrm{\scriptsize{III}}}(M)\right)\backslash
\bigcup_{x_{H}\in G}U_{x_{H},F}\}_G$$ where $G$ is a finite set; we refer to these
neighborhoods with the symbols $W_{y_{M}}$.
\item[$4^{th}$)] The point
$p_{\infty}=j_{\textrm{\scriptsize{III}}}(\omega^{\ast}\backslash\bigcup_{M\in\mathcal{M}}M^{\ast})$
which has a fundamental system of neighborhoods given by
$$\{j_{\textrm{\scriptsize{III}}}(\beta\omega)\backslash \bigcup_{y_{M}\in
K}W_{y_{M}}\}_K$$ where $K$ is a finite set.
\end{enumerate}
We could think that the points of the sets
$$
\{j_{\textrm{\scriptsize{III}}}(H^{\ast}):H\in \bigcup_{M\in
\mathcal{M}}\mathcal{H}_{M}\},\qquad \{j_{\textrm{\scriptsize{III}}}(M^{\ast}\backslash
\bigcup_{H\in {\mathcal{H}_{M}}}H^{\ast}):M\in {\mathcal{M}}\},$$$$
\{j_{\textrm{\scriptsize{III}}}(\omega^{\ast}\backslash\bigcup_{M\in\mathcal{M}}M^{\ast})\}
$$
have sequential orders respectively $1$, $2$ and $3$ with respect to the set
$j_{\textrm{\scriptsize{III}}}(\omega)$. We want to show that the point
$j_{\textrm{\scriptsize{III}}}(\omega^{\ast}\backslash\bigcup_{M\in\mathcal{M}}M^{\ast})$
does not have sequential order $3$ with respect to the set
$j_{\textrm{\scriptsize{III}}}(\omega)$. In an obvious way we point out that
$x_{\infty}\in \overline{j_{\textrm{\scriptsize{III}}}(\omega)}$: indeed there are always
infinitely many points of $\omega$ out of the union of any finite number of neighborhoods
of points of the third type; otherwise there would not be space enough for the other
elements of the MAD family $\mathcal{M}$. Now let us fix a countably infinite set
$\{M_{n}:n\in\omega\}\subseteq \mathcal{M}$ and, for every $n\in\omega$, let us fix an
infinite subset $H_{n}\in \mathcal{H}_{M_{n}}$; moreover, for every $n\in\omega$, let us
call $z_{n}$ the unique point in $K_{\textrm{\scriptsize{III}}}$ such that
$j^{-1}_{\textrm{\scriptsize{III}}}(z_{n})=H_{n}^{\ast}$. We assert that for every
$n\in\omega$ it is possible to extract a subsequence $\{m_{n_{i}}\}_{i\in\omega}$ from
$j_{\textrm{\scriptsize{III}}}(\omega)$ with $\{m_{n_{i}}\}_{i\in\omega}\rightarrow
z_{n}$: indeed for every $n\in\omega$ it is enough to put into the subsequence the image
under $j_{\textrm{\scriptsize{III}}}$ of a countably infinite number of points belonging
to $H_{n}$. Therefore for every $n\in\omega$ it turns out that $z_{n}\in
\emph{seqcl}_{1}(j_{\textrm{\scriptsize{III}}}(\omega))$. We claim that
$(z_{n})_{n\in\omega}\rightarrow p_{\infty}$: consider an open subset
$\Omega\subseteq\beta\omega$ such that $\omega^{\ast}\backslash \bigcup_{M\in
\mathcal{M}}M^{\ast}\subseteq \Omega$; we want to prove that the set $N=\{n\in\omega:
H_{n}^{\ast}\nsubseteq \Omega\}$ is finite. Towards a contradiction, suppose that $N$ is
infinite; then, in particular, it turns out that the set $\mathcal{M}'=\{M\in\mathcal{M}:
\exists H\in\mathcal{H}_{M}, H^{\ast}\nsubseteq \Omega\}$ is infinite and hence that the
set $\mathcal{M}''=\{M\in\mathcal{M}: M^{\ast}\nsubseteq \Omega\}$ is infinite. Now
consider the infinite open cover $\mathcal{A}=\{\Omega\}\cup \{M^{\ast}\cup \omega: M\in
\mathcal{M}''\}$ of $\beta\omega$; from this open cover it is not possible to extract a
finite subcover since the set $\mathcal{M}''$ is infinite. A contradiction. We can
conclude that $p_{\infty}\in \emph{seqcl}_{2}(j_{\textrm{\scriptsize{III}}}(\omega))$ and
hence that it does not have sequential order $3$ with respect to the set
$j_{\textrm{\scriptsize{III}}}(\omega)$; it follows that $K_{\textrm{\scriptsize{III}}}$
does not have sequential order $3$.

\section{BA\v{S}KIROV'S IDEA}

In this section we want to explain the general scheme of the construction
suggested by Ba\v{s}kirov's in \cite{baskirov}; the idea is to work in a
completely different way to produce compact sequential spaces of order a
successor ordinal and compact sequential spaces of order a limit ordinal. The
construction of the compact spaces of order a successor ordinal will be carried
out by transfinite induction on the order of sequentiality. When we will have
constructed compact spaces $K_{\alpha+1}$ of sequential order
$\alpha+1<\omega_{1}$ for every successor ordinal less than $\omega_1$, then it
will be easy to get compact spaces of sequential order a limit ordinal $\beta$
for every $\beta\leq\omega_{1}$. Indeed, if we want to construct a compact
sequential space of order a limit ordinal
$\beta=\sup_{\alpha+1<\beta}\{\alpha+1\}$, we can consider the disjoint sum of
the spaces $K_{\alpha+1}$ i.e.
$$Z_{\beta}=\bigoplus_{\alpha+1<\beta}
K_{\alpha+1}\ .$$ Trivially the sequential order of this space is $\beta$. It
is easy to prove that also its one-point compactification
$K_{\beta}=Z_{\beta}^{\ast}$ has also sequential order $\beta$: indeed the
point $\infty$ we have added to make compact the space has sequential order $1$
with respect to each subset
$A\subseteq K_{\beta}$ such that $\infty\in \overline{A}$.\\
 Therefore the problem reduces to construct
compact spaces whose sequential order is a successor ordinal number.
 We will construct a compact space $K_{\alpha+1}$ of
sequential order $\alpha+1$ for every successor ordinal number
$\alpha+1<\omega_{1}$; each $K_{\alpha+1}$ will be a quotient space of
$\beta\omega$, i.e. $K_{\alpha+1}=\beta\omega/\approx_{\alpha+1}$ where the
relation $\approx_{\alpha+1}$ is such that only natural numbers are one-point
elements of the quotient. For every $\alpha+1<\omega_1$ we will denote by
$j_{\alpha+1}$ the natural quotient mapping
$j_{\alpha+1}:\beta\omega\rightarrow K_{\alpha+1}$. We will prove by
transfinite induction that for each $\alpha+1<\omega_1$ the space
$K_{\alpha+1}$ satisfies the following conditions.
\begin{itemize}
  \item [S.1] The space $K_{\alpha+1}$ can be uniquely represented
in the form of $$K_{\alpha+1}=L_{0}\bigsqcup \big(\bigsqcup_{\gamma\leq\alpha}
L_{\gamma+1}\big).$$
 The points of level
$\gamma+1$ with $\gamma\in[0,\alpha]$, i.e. the points belonging
     to the set $L_{\gamma+1}$, have sequential order equal to
    $\gamma+1$ with respect to $L_{0}$, the subset consisting
     of the images of the points of $\omega$ under $j_{\alpha+1}$.
 \item [S.2] The set
$L_{\alpha+1}$ consists of only one point.
 \item [S.3] Every point in $K_{\alpha+1}$ of nonzero level has a basis formed by clopen
     subsets called elementary; moreover if $U$ is an elementary neighborhood of a point of level
     $\gamma+1$, then the relation $\approx_{\alpha+1}$ restricted to
    $\widetilde{U}=j^{-1}_{\alpha+1}(U)$ produces a compact space homeomorphic to $K_{\gamma+1}$.
  \item [S.4] For every $\gamma\leq\alpha$, if a nonconstant sequence $(x_n)_{n\in\omega}$ of points
  $x_{n}\in{L_{\gamma_{n}+1}}$, with nondecreasing levels,
   converges to a point $x\in{L_{\gamma+1}}$, then for the sequence $(\gamma_{n}+1)_{n\in\omega}$
    of ordinal numbers it holds that
    $\sup\{\gamma_{n}+1\}=\gamma$.
    \item [S.5] For every $\gamma\leq\alpha$, from every injective sequence $(x_{n})_{n\in\omega}$ of points
        $x_{n}\in{L_{\gamma_{n}+1}}$ with nondecreasing levels such that
      $\sup_{n\in\omega}\{\gamma_{n}+1\}=\gamma$, it is possible to extract
       a subsequence converging to a point of level $\gamma+1$.
  \item [S.6] If $\left\{N_{i}\right\}_{i\in\omega}$ is a countable family
       of pairwise disjoint infinite subsets $N_{i}$
     of $\omega$ and if it holds that for every $i\in\omega$ a relation of type $\beta_{i}+1$ is given on
     $\overline{N_{i}}$ in such a way that the sequence of ordinals $(\beta_{i}+1)_{i\in\omega}$ is not
     decreasing and $\sup\left\{\beta_{i}+1\right\}=\alpha$, then it is possible
     to extend the relation obtained on $\bigcup_{i=1}^{\infty}\overline{N_{i}}$ to a relation of $\beta\omega$ of type
     $\alpha+1$.
\end{itemize}

From the first three conditions we trivially deduce other two
properties.
\begin{itemize}
\item [S.7] If $U$ is an elementary neighborhood of a point $x$ of
level $\gamma+1$ in $K_{\alpha+1}$, then its level in
$U=\widetilde{U}/(\approx_{\alpha+1}\ |_{\widetilde{U}})$ is equal
to $\gamma+1$. \item [S.8] If $U$ is an elementary neighborhood of
a point $x$ of level $\gamma+1$ in $K_{\alpha+1}$, then
 $U\backslash\left\{x\right\}\subseteq\bigcup_{\gamma'<\gamma}L_{\gamma'+1}$.
\end{itemize}
\vspace{0,2cm}

The compact sequential spaces $K_{1}$ and $K_{2}$ will be taken as
bases of the recursion. Let us begin to check that properties S.1
to S.6 hold for these spaces.

\begin{center}\textbf{Check of the properties of $K_{1}$}
\end{center}
\begin{enumerate}
\item [S.1] The space $K_{1}=\omega\cup\left\{P\right\}$ can be
uniquely represented in the form of $K_{1}=L_{0}\bigsqcup L_{1}$; we denote by $L_{0}$
the one-point elements of the quotient that are images of the points of $\omega$ under
$j_1$, while $L_{1}$ consists of only one point that is the image of $\omega^{\ast}$
under the natural quotient mapping. The unique point of $L_{1}$ has sequential order $1$
with respect to $L_{0}$.
 \item [S.2] The set $L_{1}$ consists of the unique limit point of
 the sequence.
  \item [S.3]  The unique point of $L_{1}$ has a basis formed by the
  clopen elementary
  neighborhoods $U_{P,F}$: the space obtained by restricting the relation $\approx_{1}$ to
  $\widetilde{U}_{P,F}=j_{1}^{-1}(U_{P,F})$
  is homeomorphic to $K_{1}$.

   \item [S.4] For $K_{1}$ we do not have anything to prove.
   \item[S.5] Obvious.
    \item [S.6] It is necessary to prove S.6 starting from the space $K_{2}$.
\end{enumerate}
\vspace{0,2cm}

\begin{center}\textbf{Check of the properties of $K_{2}$}
\end{center}

\begin{enumerate}
\item [S.1] The space $K_{2}$ can be uniquely represented in the
form of $$K_{2}=L_{0}\bigsqcup L_{1}\bigsqcup L_{2}$$ where we
denote by $L_{0}$ the one-point elements of the quotient that are
images of the points of $\omega$ under $j_2$; the quotient mapping
$j_{2}$ collapses every $M^{\ast}$ with $M\in \mathcal{M}$ to a
single point (and $L_{1}$ consists of these points), while it
collapses
$\omega^{\ast}\backslash\bigcup_{M\in{\mathcal{M}}}M^{\ast}$ to
another single point that gives $L_{2}$.\\
 The points of $L_{1}$
have sequential order $1$ with respect to $L_{0}$, while the unique point of $L_{2}$ has
sequential order $2$ with respect to $L_{0}$.
\item [S.2] The set $L_{2}$ consists of the unique point
$$x_{\infty}=j_2(\omega^{\ast}\backslash\bigcup_{M\in{\mathcal{M}}}M^{\ast}).$$
 \item [S.3] Every point $y\in L_1$ has a basis formed by the
 clopen elementary neighborhoods $U_{y,F}$ and the space obtained by restricting the relation
  $\approx_{2}$ to $\widetilde{U}_{y,F}=j_{2}^{-1}(U_{y,F})$ is homeomorphic
   to the compact sequential space $K_{1}$.
     The point $x_{\infty}\in L_2$ has a basis given by the
     clopen elementary
  neighborhoods
   $K_{2}\backslash \bigcup_{x\in G} U_{x}$; the space obtained by restricting the relation
  $\approx_{2}$ to $j_{2}^{-1}(K_{2}\backslash \bigcup_{x\in G} U_{x})$
   is a compact space homeomorphic to $K_{2}$.
  \item [S.4] For $K_{2}$ Property S.4 is obvious.
  \item[S.5] From every noncostant sequence $(x_{n})_{n\in\omega}$ of points with nondecreasing levels such that
      $\sup\{l(x_{n})\}=0$ it is possible to extract
       a subsequence converging to a point of level $1$: indeed, since the family $\mathcal{M}$ is
       MAD,
       the sequence $(x_{n})_{n\in\omega}$ has infinite intersection with at least one element $M_{1}$ of the family
       $\mathcal{M}$ and hence we can extract a
       subsequence converging to the point $j_{2}(M_{1}^{\ast})$ of level $1$.
From every noncostant sequence
      $(x_{n})_{n\in\omega}$ of points with nondecreasing levels such that
      $\sup\{l(x_{n})\}=1$ it is possible to extract
       a subsequence converging to $x_{\infty}$:
      indeed the points of the sequence are eventually in
      every neighborhood of $x_{\infty}$.
  \item [S.6] If $\left\{N_{i}\right\}_{i\in\omega}$ is a countable family
       of pairwise disjoint infinite subsets $N_{i}\subset\omega$ and if it holds that for every $i\in\omega$
        a relation of type $\beta_{i}+1$ is given on
     $\overline{N_{i}}$ in such a way that the sequence of ordinals $\beta_{i}+1$ is
     nondecreasing and that $\sup\left\{\beta_{i}+1\right\}=1$ (and hence in such a way that $\beta_{i}+1=1$
for every $i\in\omega$),
     then we can extend the so obtained relation on $\bigcup_{i=1}^{\infty}\overline{N_{i}}$
      to a relation on $\beta\omega$ of type
     $2$: indeed it is enough to complete the almost disjoint family $\{N_{i}\}_{i\in\omega}$
      to a MAD family and then put into relation the elements of $\beta\omega$
       in the way we have already seen when we constructed $K_2$.
     \end{enumerate}

Now we are sure we can take the compact sequential spaces $K_1$ and $K_2$ as
bases of the induction. Moreover we will assume that for all $\beta+1<\alpha+1$
the compact sequential spaces $K_{\beta+1}$ (in which properties S.1 to S.6
hold) have been constructed; then we will able to construct the compact space
$K_{\alpha+1}$ with sequential order $\alpha+1$ satisfying conditions S.1 to
S.6.
\section{SOME PROPAEDEUTIC LEMMAS} Before showing the construction in
all details, let us prove some useful very technical lemmas.

\begin{lemma}\label{lemma0}
The intersection of any countable family of open subsets of
$\omega^{\ast}$ is either empty or contains a non-empty open
subset.
\end{lemma}
\iniziodim Let $\{A_{i}\}_{i\in\omega}$ be a countable family of open subsets
of $\omega^{\ast}$ whose intersection contains a point $\mathcal{U}$. For every
$i\in\omega$ there exists a subset $N_{i}\subset\omega$ such that
$$\mathcal{U}\in N_{i}^{\ast}\subset A_i .$$
The intersection of any finite collection of the sets $N_{i}^{\ast}$ is not empty and
open and hence the intersection of any finite collection of the sets $N_{i}$ is infinite.
Then there exists an increasing sequence of integers $n_i$ such that $n_i\in N_{1}\cap
N_{2}\cap\ldots\cap N_{i}$. Let us set $N=\{n_i:i\in\omega\}$; since $N\backslash N_{i}$
is finite for each $i\in\omega$, it holds that $N^{\ast}\subset N_{i}^{\ast}$ for each
$i\in\omega$ and then we can conclude that $N^{\ast}\subset \bigcap_{i} A_{i}$ where
$N^{\ast}\neq\emptyset$ as $|N|=\omega$.
 \finedim

\begin{lemma}\label{lemma01} Let $\{N_{i}^{\ast}\}_{i\in\omega}$ be a countably infinite family of
 clopen subsets of $\omega^{\ast}$; let us
suppose that
$\omega^{\ast}\backslash\overline{\bigcup_{i\leq\bar\imath}N_{i}^{\ast}}\neq\emptyset$
for every $\bar\imath\in\omega$. Then there exists
$\Delta\subset\omega$ such that $|\Delta|=\omega$ and
$\Delta^{\ast}\cap \bigcup_{i\in\omega}N_{i}^{\ast}=\emptyset$.
\end{lemma}
\iniziodim Let us set $\Delta_{1}=N_{1}^{\ast},\Delta_{2}=N_{1}^{\ast}\cup
N_{2}^{\ast},\ldots,\Delta_{\bar\imath}=N_{1}^{\ast}\cup
N_{2}^{\ast}\cup\ldots\cup N_{\bar\imath}^{\ast}$ and so on; it is clear that
for every $\bar\imath\in\omega$, $\Delta_{\bar\imath}$ is a clopen subset of
$\omega^{\ast}$. Notice that $\omega^{\ast}\backslash \Delta_{1}\supseteq
\omega^{\ast}\backslash \Delta_{2}\supseteq \ldots \supseteq
\omega^{\ast}\backslash \Delta_{\bar\imath}\supseteq \omega^{\ast}\backslash
\Delta_{\bar\imath+1}\supseteq\ldots$; let us set
$C_{i}=\omega^{\ast}\backslash \Delta_{i}$ for every $i\in\omega$ and
$\mathcal{E}=\{C_{i}:i\in\omega\}$. The family $\mathcal{E}$ satisfies the
finite intersection property i.e. the intersection of a finite number of
elements of $\mathcal{E}$ is not empty: indeed the set of
 the indices of this finite number of elements has a maximum
$\overline{n}$ and their intersection is equal to
$\omega^{\ast}\backslash\bigcup_{j\leq\overline{n}}N_{j}^{\ast}\supseteq
\omega^{\ast}\backslash\overline{\bigcup_{j\leq\overline{n}}N_{j}^{\ast}}$ which is a
non-empty subset by hypothesis. Since $\omega^{\ast}$ is compact, it follows that
$\bigcap \mathcal{E}=\bigcap_{i\in\omega}C_{i}\neq\emptyset$. Therefore it turns out that
\begin{displaymath}\emptyset\neq\bigcap_{i\in\omega}C_{i}=\bigcap_{i\in\omega}\Delta_{i}^{C}=
[\bigcup_{i\in\omega}(\cup_{j\leq
i}N_{j}^{\ast})]^{C}=(\bigcup_{i\in\omega}N_{i}^{\ast})^{C}=\omega^{\ast}\backslash
\bigcup_{i\in\omega}N_{i}^{\ast}.
\end{displaymath}
We can conclude that the family $\{\Delta_{i}^{C}\}_{i\in\omega}$
consisting of open subsets of $\omega^{\ast}$ has non-empty
intersection and hence, by Lemma \ref{lemma0}, this intersection
contains an open subset $A\subseteq \omega^{\ast}\backslash
\bigcup_{i\in\omega}N_{i}^{\ast}$. Then there exists
$\Delta\subset\omega$ with $|\Delta|=\omega$ such that
$\Delta^{\ast}\subseteq\omega^{\ast}\backslash
\bigcup_{i\in\omega}N_{i}^{\ast}$ and hence such that
$\Delta^{\ast}\cap \bigcup_{i\in\omega}N_{i}^{\ast}=\emptyset$.
 \finedim

\begin{lemma}\label{lemma2} Let
$\mathcal{P}=\mathcal{Q}\cup\mathcal{R}$ be a family of infinite subsets of $\omega$ such
that
\begin{itemize}
 \item[-] $\mathcal{Q}$ is an almost disjoint family;
\item[-] $|\mathcal{Q}|\leq\omega$ and $|\mathcal{R}|\leq\omega$;
\item[-] for every element $Q_{i}\in\mathcal{Q}$ and every element
$R_{n}\in\mathcal{R}$ it turns out that $|Q_{i}\cap R_{n}|<\omega$.
\end{itemize}
Then there exists $L\in[\omega]^{\omega}$ such that
$L^{\ast}\supseteq\bigcup_{Q_{i}\in\mathcal{Q}}Q_{i}^{\ast}$ and $L^{\ast}\cap
R_{n}^{\ast}=\emptyset$ for every $R_{n}\in\mathcal{R}$.
\end{lemma}

\iniziodim Let us set $\mathcal{Q}=\{Q_{i}: i\in\omega\}$ with $|Q_{i}\cap
Q_{j}|<\omega$ for $i\neq j$. Obviously we can suppose
$\mathcal{R}\neq\emptyset$ and then we can write $\mathcal{R}=\{R_{n}:
n\in\omega\}$ with $n\mapsto R_{n}$ not necessarily injective. If
$|\mathcal{Q}|<\omega$, then we set $L=\bigcup Q_{i}$. If instead
$|\mathcal{Q}|=\omega$, we set
$L=\bigcup_{n\in\omega}(Q_{n}\backslash\cup_{n'<n}R_{n'})$. For every
$\overline{n}\in\omega, R_{\overline{n}}$ intersects $L$ only in those points
in which $R_{\overline{n}}$, in case, intersects the subsets $Q_{n}$ with
$n=0,\ldots,\overline{n}$; these points are in a finite number. Furthermore for
every $\overline{n}\in\omega$, $Q_{\overline{n}}\backslash L$ consists of a
finite number of points and exactly of those in which $Q_{\overline{n}}$
intersects $R_{n}$ with $n<\overline{n}$ (and these again are certainly in a
finite number); therefore $Q_{\overline{n}}^{\ast}\subseteq L^{\ast}$ for every
$\overline{n}\in \omega$ and hence $\bigcup_{n}Q_{n}^{\ast}\subseteq
L^{\ast}$.\finedim

In the following lemma we will take into account a countable
family of infinite pairwise disjoint subsets of $\omega$,
$\{\tilde{N}_{i}\}_{i\in\omega}$ and a relation $\approx$ on
$U=\bigsqcup \tilde{N}_{i}^{\ast}\subset\omega^{\ast}$. Let us set
$H=U/\approx$ and let $j$ be the quotient mapping $j: U\rightarrow
H$. We assume that the subsets $\tilde{N}_{i}^{\ast}$ are
\emph{distinguished relative to $\approx$}, i.e. that
$j^{-1}(j(\tilde{N}_{i}^{\ast}))=\tilde{N}_{i}^{\ast}$ for every
$i\in\omega$. Now let us prove the lemma.

\begin{lemma}\label{lemma1} Let $\{\tilde{N}_{i}\}_{i\in\omega}$ be a countable family
 of infinite pairwise disjoint subsets
$\tilde{N}_{i}\subset\omega$ and let $\approx$ be a relation on $U=\bigsqcup
\tilde{N}_{i}^{\ast}\subset\omega^{\ast}$ where the subsets
$\tilde{N}_{i}^{\ast}$ are distinguished relative to $\approx$ and the spaces
$\tilde{N}_{i}^{\ast}/\approx$ are zero-dimensional compact spaces. Let us
suppose that the set $B=\{x_{n}:n\in\omega\}$ is devoid of any accumulation
point in $H$ and that for every $x_{n}$ there exists an index $i_{n}$ and a
clopen neighborhood $U(x_{n})$ such that $U(x_{n})\subseteq
\tilde{N}_{i_{n}}^{\ast}/\approx$. Moreover let us suppose that
$\bigcup_{n\in\omega} U(x_{n})\neq H$. Then
\begin{itemize}
\item[i)] there exist pairwise disjoint clopen subsets $U_{n}$
with $n\in\omega$ such that $x_{n}\in U_{n}$ for every
$n\in\omega$; \item[ii)] there exists $\tilde{N'}\subset\omega$
such that
$$(\tilde{N'})^{\ast}\cap \bigcup_{i\in\omega} \tilde{N}_{i}^{\ast}= j^{-1}
(\bigsqcup_{n\in\omega} U_{n})=\bigsqcup_{n\in\omega} E_{n}^{\ast}.$$
\end{itemize}
 \end{lemma}
\iniziodim Since the subsets $\{\tilde{N}^{\ast}_{i}\}_{i}$ are disjoint and
distinguished relative to $\approx$
 it holds that $$(\tilde{N}_{i}^{\ast}/\approx)\cap
(\tilde{N}_{j}^{\ast}/\approx)=\emptyset$$ for every $i,j\in\omega$ with $i\neq
j$. Moreover $\tilde{N}_{i}^{\ast}/\approx$ is open and closed in $H$ for every
$i\in\omega$, since $j^{-1}(\tilde{N}_{i}^{\ast}/\approx)=\tilde{N}_{i}^{\ast}$
is open and closed in $U$. We need to remark that, for every $i\in\omega$,
$\tilde{N}_{i}^{\ast}/\approx$ intersects only a finite number of the
neighborhoods $\{U(x_n)\}_{n}$ because of the hypothesis that the subset B is
devoid of any accumulation point in $H$.\\
 By transfinite induction we are going to construct
the subsets
$U_{n}$ with $n\in\omega$ such that $x_{n}\in U_{n}$ for every $n\in\omega$.\\
Let us consider the point $x_{1}$ and the clopen subset $U(x_{1})\subseteq
\tilde{N}_{i_{1}}^{\ast}/\approx$; since $B$ is devoid of any accumulation
point in $H$, there exists an open neighborhood $A_{1}\subseteq H$ of $x_{1}$
such that $A_{1}\cap B=\{x_{1}\}$. Now $x_{1}\in [(A_{1}\cap
\tilde{N}_{i_{1}}^{\ast}/\approx)\cap U(x_{1})]=D_{1}$: this subset is open in
$\tilde{N}_{i_{1}}^{\ast}/\approx$ and hence, since
$\tilde{N}_{i_{1}}^{\ast}/\approx$ is zero-dimensional, there exists a clopen
subset $U_{1}\subseteq D_{1}$ of $\tilde{N}_{i_{1}}^{\ast}/\approx$ with
$x_{1}\in U_{1}$; trivially $U_{1}$ is clopen also in $H$. Now $H\backslash
U_{1}$ is an open subset of $H$ and it contains $B\backslash \{x_{1}\}$ which
is devoid of any accumulation point in $H\backslash U_{1}$; thus there exists
an open neighborhood $A_{2}\subseteq H$ of $x_{2}$ such that $A_{2}\cap
B=\{x_{2}\}$. Hence $x_{2}\in [(H\backslash U_{1})\cap A_{2}\cap
(\tilde{N}_{i_{2}}^{\ast}/\approx)\cap U(x_{2})]=D_{2}$: this is an open subset
of $\tilde{N}_{i_{2}}^{\ast}/\approx$ and then, since
$\tilde{N}_{i_{2}}^{\ast}/\approx$ is zero-dimensional, there exists a clopen
subset $U_{2}\subseteq D_{2}$ of $\tilde{N}_{i_{2}}^{\ast}/\approx$ with
$x_{2}\in U_{2}$; trivially $U_{2}$ is clopen also in $H$ and moreover it
results that $U_{1}\cap U_{2}=\emptyset$. Notice that $H\backslash (U_{1}\sqcup
U_{2})$ is a non-empty clopen subset of $H$ and that $B\backslash
\{x_{1},x_{2}\}\subseteq H\backslash
(U_{1}\sqcup U_{2})$.\\
 Let us suppose
that for every $n\leq \overline{n}$ there exists a clopen subset
$U_{n}\subseteq H$ with $x_n\in U_{n}$ and $U_{n}\subseteq U(x_{n})$ and that
$U_{n}\cap U_{n'}=\emptyset$ for every $n',n\leq\overline{n}$; moreover suppose
that for every $n\leq\overline{n}$ it holds that $B_{n}=B\backslash
\{x_{1},\ldots,x_{n}\}\subseteq H\backslash \bigsqcup_{j\leq n}U_{j}$. Let us
prove that these properties hold also for $\overline{n}+1$. By inductive
hypothesis $x_{\overline{n}+1}\in(H\backslash \bigsqcup_{j\leq
\overline{n}}U_{j})$ where $H\backslash \bigsqcup_{j\leq \overline{n}}U_{j}$ is
open, since the finite union $\bigsqcup_{j\leq \overline{n}}U_{j}$ is clopen;
moreover there exists an open neighborhood $A_{\overline{n}+1}\subseteq H$ of
$x_{\overline{n}+1}$ such that $A_{\overline{n}+1}\cap
B=\{x_{\overline{n}+1}\}$. It turns out that $x_{\overline{n}+1}\in
[(H\backslash \bigsqcup_{j\leq n}U_{j})\cap A_{\overline{n}+1}\cap
\tilde{N}_{i_{\overline{n}}+1}^{\ast}/\approx] \cap
U(x_{\overline{n}+1})=D_{\overline{n}+1}$ and that $D_{\overline{n}+1}$ is open
in $\tilde{N}_{i_{\overline{n}}+1}^{\ast}/\approx$. Now, since
$\tilde{N}_{i_{\overline{n}}+1}^{\ast}/\approx$ is zero-dimensional, there
exists a clopen subset $U_{\overline{n}+1}\subseteq D_{\overline{n}+1}$ of
$\tilde{N}_{i_{\overline{n}}+1}^{\ast}/\approx$ with $x_{\overline{n}+1}\in
U_{\overline{n}+1}$: we trivially remark that $U_{\overline{n}+1}$ is also
clopen in $H$, that $U_{\overline{n}+1}\cap U_{n'}=\emptyset$ for every
$n'<\overline{n}+1$ and that $H\backslash \bigsqcup_{j\leq
\overline{n}+1}U_{j}\supseteq B\backslash
\{x_{1},\ldots,x_{\overline{n}+1}\}$.\\
 Therefore $U_{n}\subseteq
\tilde{N}_{i_{n}}^{\ast}/\approx$ is a clopen neighborhood of
$x_n$ in $H$ for every $n\in\omega$ and $j^{-1}(U_{n})$ is a
clopen subset of $\tilde{N}_{i_{n}}^{\ast}$; then for every
$n\in\omega$ it holds that $j^{-1}(U_{n})=E_{n}^{\ast}$ where
$E_{n}$ is an infinite subset of $\omega$.
 Moreover we can assert that $\bigsqcup_{n\in\omega}
U_{n}\neq H$ since
$\bigsqcup_{n\in\omega}U_{n}\subseteq\bigcup_{n\in\omega}U(x_{n})$. Now we want
to prove that $\bigsqcup_{n\in\omega} U_{n}$ is clopen in $U/\approx$;
trivially $\bigsqcup_{n\in\omega}U_{n}$ is open and now we show that it is also
closed. If we take a point $z\in H\backslash \bigsqcup_{n\in\omega} U_{n}$
there exists an index $i_{z}$ such that
$z\in(\tilde{N}_{i_{z}}^{\ast}/\approx)\backslash \bigsqcup_{n\in\omega}
U_{n}$.
 Since $\tilde{N}_{i_{z}}^{\ast}/\approx$ intersects only a finite
number of the clopen subsets we have just constructed (we denote these clopen
subsets by $U_{j_{1}},\ldots,U_{j_{\overline{n}}})$, then
$\bigsqcup_{i=1}^{\overline{n}}U_{j_{i}}\cap
(\tilde{N}_{i_{z}}^{\ast}/\approx)$ is closed in
$\tilde{N}_{i_{z}}^{\ast}/\approx$ since
$\bigsqcup_{i=1}^{\overline{n}}U_{j_{i}}$ is closed in $H$; the subset
$\tilde{N}_{i_{z}}^{\ast}/\approx\backslash(\bigsqcup_{i=1}^{\overline{n}}U_{j_{i}})$
is an open subset to which $z$ belongs and hence there exists an open
neighborhood of $z$ in $\tilde{N}_{i_{z}}^{\ast}/\approx$ (and then in $H$)
disjoint from $\bigsqcup U_{n}$. Finally we can conclude that $j^{-1}(\bigsqcup
U_{n})=\bigsqcup j^{-1}(U_{n})$ is clopen in $U$ and then there exists a clopen
subset $(\tilde{N}')^{\ast}$ of $\omega^{\ast}$ with
$\tilde{N}'\subseteq\omega$ and $|\tilde{N}'|=\omega$ such that
\begin{displaymath}(\tilde{N}')^{\ast}\cap \bigcup \tilde{N}_{i}^{\ast}= \bigsqcup
j^{-1}(U_{n})=\bigsqcup E_{n}^{\ast}.
 \end{displaymath}
\finedim

\begin{remark}\label{ossdiscreto} We remark that Lemma \ref{lemma1}
 still holds when we consider a family
 $$\left\{\tilde{N}_{\gamma}:\gamma\in\omega_{1}\right\}$$
  of infinite subsets $\tilde{N}_{\gamma}\subset\omega$ keeping all the other
  hypotheses.
\end{remark}

\section{Construction of a Ba\v{s}kirov's space of order an arbitrary successor ordinal}

Finally we will show how to construct the space $K_{\alpha+1}$
 by
assuming that all compact sequential spaces $K_{\beta+1}$ (of
sequential order $\beta+1$) with $\beta+1<\alpha+1$ have been
constructed and that properties S.1 to S.6 hold in each of these
space; moreover we will check that properties S.1 to S.6 hold in
$K_{\alpha+1}$ too. We will carry out the construction when
$\alpha$ is a successor ordinal, but we will remark from time to
time what it is necessary to change if
we have to work in the case in which $\alpha$ is a limit ordinal).\\
It will be very important to take the set $\Gamma$ into account: it is the set
of all the families $\mathcal{C}_{\xi}$ whose elements are countable pairwise
disjoint clopen subsets of $\omega^{\ast}$; under the Continuum Hypothesis, we
can write $\Gamma$ as
\begin{equation}\label{gammabaski}
 \Gamma=\left\{\mathcal{C}_{\xi}:\omega\leq
\xi <\omega_{1}\right\}.
\end{equation}
 Roughly speaking, our type of construction ensures that,
  by a number of steps of cardinality equal to the cardinality of
$\Gamma$, we are able to exhaust the whole $\Gamma$; moreover at each stage
$\alpha<\omega_1$ of the inductive construction, it will be essential the fact that
$\alpha$ is a countable ordinal in order to guarantee that we can continue the process
and hence it is
crucial that we can enumerate $\Gamma$ as in (\ref{gammabaski}).\\
 Let us begin the construction.
 Let $\left\{N_{i}\right\}_{i\in\omega}$ be a family of pairwise
disjoint infinite subsets $N_{i}\subset\omega$. For every $i\in\omega$ the
closures of $N_i$ in $\beta\omega$, namely $\overline{N_{i}}$, is a clopen
subset of $\beta\omega$ which is homeomorphic to it; for every $i\in\omega$ let
us set a decomposition of type $\beta_{i}+1$ on $\overline{N_{i}}$ taking care
that the sequence of ordinals $S=(\beta_{i}+1)_{i\in\omega}$ is nondecreasing
and such that $\sup \{\beta_{i}+1\}=\alpha$. Notice that it is possible to
extract an injective subsequence $S'=(\beta_{i_n}+1)_{n\in\omega}\subseteq S$
in such a way that the sequence $S'$ converges upwards to $\alpha$; since
$\alpha$ is a successor ordinal, this mean that there are infinite $n\in\omega$
such that the decomposition set on $\overline{N_{i_n}}$ is a relation of type
$\alpha$.\footnote{If $\alpha$ is a limit ordinal we will have to set
decompositions of type $\beta_{i}+1$ on $\overline{N_{i}}$ in such a way that
the sequence $(\beta_{i}+1)_{i\in\omega}$ is nondecreasing and $\sup
\{\beta_{i}+1\}=\alpha$. Also in this case it is possible to extract an
injective subsequence $S'=(\beta_{i_n}+1)_{n\in\omega}\subseteq S$ in such a
way that the sequence $S'$ converges upwards to
$\alpha$.}\\
For every $i\in\omega$, let $j_{\beta_{i}+1}:\overline{N_{i}}
\rightarrow K_{\beta_{i}+1}$ be the quotient mapping. Let us check
that the following properties hold for every
$\bar\imath\in\omega$.
\begin{enumerate}
    \item [T.1]
    $\overline{N_{\bar\imath}^{\ast}\backslash\bigcup_{i'<\ \bar\imath}
    N_{i'}^{\ast}}\neq \emptyset$: indeed the subsets $N_{i}$
     are pairwise disjoint and, in particular, the subsets $N_{i}$
     with $i=1,\ldots,\bar\imath$ are pairwise disjoint; then
     the corresponding $N_{i}^{\ast}$ are pairwise disjoint.
      Hence it holds that $N_{\bar\imath}^{\ast}\backslash (\bigcup_{i'<\ \bar\imath}
    N_{i'}^{\ast})\neq \emptyset$ and obviously
     $\overline{N_{\bar\imath}^{\ast}\backslash\bigcup_{i'<\ \bar\imath}
    N_{i'}^{\ast}}\neq \emptyset$.

    \item [T.2] $\overline{\bigcup_{i'\leq \ \bar\imath}N^{\ast}_{i'}}\neq \omega^{\ast}$:
     indeed we have already pointed out that the subsets $N_{i}^{\ast}$
      are pairwise disjoint; then $N_{\bar\imath+1}^{\ast}$
      is a non-empty subset disjoint from all the subsets $N_{i}^{\ast}$ with $i\leq
     \bar\imath$.
Moreover for every $i\leq\bar\imath$ the subset $N_{i}^{\ast}$ is
closed and the finite union $\bigcup_{i\leq\
\bar\imath}N_{i}^{\ast}$ is closed too; thus it follows that
$\omega^{\ast}\backslash \overline{\bigcup_{i\leq\
\bar\imath}N_{i}^{\ast}}=\omega^{\ast}\backslash \bigcup_{\leq\
\bar\imath}N_{i}^{\ast}\supset N_{\bar\imath+1}^{\ast}$.

\item [T.3] For every $i'<\bar\imath$ it holds that
$\overline{N_{i'}}\cap\overline{N_{\bar\imath}}=\emptyset$.
    \item [T.4] Property T.4 takes into account the families
    $\mathcal{C}_{\xi}$ with $\xi\leq\bar\imath$ but the families
    $\mathcal{C}_{\xi}$ have indices from $\omega$ to
    $\omega_{1}$ not included and so, at the moment, we do not have to consider
    this property.
\end{enumerate}

 In view
of T.3 and the relations set on each $\overline{N_{i}}$ with
$i\in\omega$, we have defined a relation $Q_{\omega}$ on
$U_{\omega}=\bigcup_{i\in\omega}N_{i}^{\ast}$. Let
$j^{\omega}_{\alpha+1}$ be the quotient mapping
$j^{\omega}_{\alpha+1}:U_{\omega}\rightarrow
U_{\omega}/Q_{\omega}$.\\ We say that
$\mathcal{C}_{\xi}\in{\Gamma}$ is an $\omega-$family if
$\mathcal{C}_{\xi}$ consists of elements that can be decomposed
into two subfamilies $\mathcal{L}_{0}$ and $\mathcal{L}_{1}$
satisfying the following conditions.
\begin{enumerate}
\item [U.1] $\bigcup \mathcal{L}_{0}\cap U_{\omega}=\emptyset$.
 \item [U.2] For every
$c\in{\mathcal{L}_{1}}$ there exists $i<\omega$, a point
$x_{c}\in{\overline{N_{i}}/\approx_{\beta_{i}+1}}$ of level
 $\gamma_{c}+1$ and an elementary neighborhood $U_{c}$ of $x_{c}$
 such that
 $c=\widetilde{U_{c}}\cap\omega^{\ast}$ where $\widetilde{U_{c}}=j_{\beta_{i}+1}^{-1}(U_{c})$.
 \item [U.3] The set $\left\{x_{c}:
c\in{\mathcal{L}_{1}}\right\}$ has no accumulation points in
$U_{\omega}/Q_{\omega}$. \item [U.4] It holds that
$\sup\left\{\gamma_{c}+1:c\in{\mathcal{L}_{1}}\right\}<\alpha$.
\end{enumerate}
Let us rewrite these properties in order to make clear the new notion.
\begin{enumerate}
\item[i)] $\mathcal{L}_{0}$ consists of elements $C_{n}^{\ast}$
where the subsets $C_{n}\subset\omega$ are transversal to the
subsets $N_{i}$, i.e. every $C_{n}\subset\omega$ intersects every
$N_{i}$ in a finite number of points (in this way we are
respecting U.1); \item[ii)] $\mathcal{L}_{1}$ consists of elements
$C_{m}^{\ast}$ where for every $m$ there exist $i\in\omega$, a
point $x_{m}\in \overline{N_{i}}/\approx_{\beta_{i}+1}$ of level
$l(x_{m})<\alpha$ and an elementary neighborhood $U(x_{m})$ such
that $C_{m}^{\ast}=j_{\beta_{i}+1}^{-1}U(x_{m})\cap\omega^{\ast}$.
A further necessary requirement is that the set $\{x_{m}\}$ is
devoid of any accumulation point in $U_{\omega}/Q_{\omega}$ and
that $\sup\{l(x_{m})\}<\alpha$. We want to point out that by
$l(x_{j})$ we mean a successor ordinal.
 (In this way we are respecting U.2-U.3-U.4). \end{enumerate}
Notice that it is possible to find an $\omega$-family: for
example, we can use Lemma \ref{lemma01} since the subsets
$N_{i}^{\ast}$ comply with the hypotheses; in this way we find an
infinite subset $\Delta_{\omega}\subset\omega$ such that
$\Delta_{\omega}$ intersects every $N_{i}$ in a finite number of
points. We can decompose this infinite set in an infinite number
of infinite subsets $T_n\subset\omega$ that again intersect every
$N_{i}$ in a finite number of points; we set
$\mathcal{L}_{0}=\{T_n^{\ast}:n\in\omega\}$. It is clear that
$\mathcal{L}=\mathcal{L}_{0}$ is an $\omega$-family.\\
 Among all the $\omega$-families let us take
the one with the minimum index $\overline{\omega}$; we write it as
$\mathcal{C}_{\overline{\omega}}=\mathcal{L}_{0}\sqcup\mathcal{L}_{1}$ with
$\mathcal{L}_{0}=\{N_{\overline{\omega},n}^{\ast}:n\in J_0\}$,
$\mathcal{L}_{1}=\{N_{\overline{\omega},n}^{\ast}:n\in J_1\}$ and $J_0\cap
J_1=\emptyset$. Of course, by construction, the $\omega$-family $C_{\overline{\omega}}$
complies with the following properties.
\begin{enumerate}
\item [U.1] $\bigcup \mathcal{L}_{0}\cap U_{\omega}=\emptyset$.
 \item [U.2] For every $N_{\overline{\omega},n}^{\ast}$
  with $n\in J_{1}$ there exist $i_{n}\in\omega$,
 a point $x_{n}
  \in \overline{N_{i_{n}}}/\approx_{\beta_{i_{n}}+1}$ of level $l(x_{n})<\alpha$
   and an elementary neighborhood
   $U(x_{n})$ such that
   $N_{\overline{\omega},n}^{\ast}=j_{\beta_{i_{n}}+1}^{-1}(U(x_{n}))\cap\omega^{\ast}$.
 \item [U.3] The set $\{x_{n}:n\in J_1\}$ has no accumulation point in $U_{\omega}/Q_{\omega}$. \item [U.4]
It holds that $\sup\{l(x_{n}):n\in J_{1}\}=\beta_{\omega}<\alpha$
with $\beta_{\omega}$ that can take up value from $1$ to $\alpha$
not included. Without loss of generality we can always assume that
the levels of the points are ordered in a nondecreasing way.
\end{enumerate}
For every $n\in J_{1}$ it turns out that
$\hat{U}(x_{n})=U(x_{n})\backslash\omega$ is a clopen neighborhood
of $x_{n}$ in $N_{i_{n}}^{\ast}/\approx_{\beta_{i_{n}}+1}$. We can
apply Lemma \ref{lemma1} since the family $\{N_{i}:i\in\omega\}$,
the points $x_{n}$ with $n\in J_{1}$ and the relation $Q_{\omega}$
defined on $U_{\omega}=\bigsqcup N_{i}^{\ast}$ satisfy the
hypotheses. We remark that $\bigcup\hat{U}(x_{n})\neq
U_{\omega}/Q_{\omega}$ since in $U_{\omega}/Q_{\omega}$ there are
points of level $\alpha$
 which $\bigcup\hat{U}(x_{n})$ does not cover.\footnote{Notice that $\bigcup\hat{U}(x_{n})\neq
U_{\omega}/Q_{\omega}$
  also in the case in which $\alpha$ is a limit ordinal:
 indeed at the beginning of the construction we put decompositions of type $\beta_{i}+1$ on
 the subsets $\overline{N}_{i}$ in such a way that
 $\sup\{\beta_{i}+1\}=\alpha$; hence in $U_{\omega}/Q_{\omega}$ there certainly
 exists a point of level $\beta_{\omega}+1<\alpha$
   that $\bigcup\hat{U}(x_{n})$ does not cover.}
 Therefore it is possible to find pairwise elementary neighborhoods $U_{n}$
with $x_{n}\in U_{n}$ and a subset $N_{\omega}'\subset \omega$ such that
$$(N_{\omega}')^{\ast}\cap U_{\omega}=\bigsqcup_{n\in J_1}
(j_{\alpha+1}^{\omega})^{-1}(U_{n})=\bigsqcup_{n\in J_1} E_{n}^{\ast}.$$ Let us define
$\mathcal{C}'=\mathcal{L}_{0}\cup \{(j_{\alpha+1}^{\omega})^{-1}(U_{n}):n\in J_{1}\}$.\\
Now if we set $\mathcal{Q}=\{N_{\overline{\omega},n}:n\in J_0\}$ and
$\mathcal{R}=\{N_i:i\in\omega\}$, then $\mathcal{P}=\mathcal{Q}\cup\mathcal{R}$ is a
family of subsets with the following properties:
\begin{itemize}
 \item[-] $\mathcal{Q}$ is an almost disjoint family;
\item[-] $|\mathcal{Q}|\leq\omega$ and $|\mathcal{R}|\leq\omega$;
\item[-] for every $N_{\overline{\omega},n}\in \mathcal{Q}$ and
every $N_{i}\in\mathcal{R}$ it holds that $|N_{\overline{\omega},n}\cap N_{i}|<\omega.$
\end{itemize}
Therefore, by Lemma \ref{lemma2}, there exists a subset
$N_{\omega}''\in[\omega]^{\omega}$ such that $$\bigcup_{n\in
J_{0}}N_{\overline{\omega},n}^{\ast}\subseteq (N_{\omega}'')^{\ast} \quad
\textmd{and}\quad (N_{\omega}'')^{\ast}\cap N_{i}^{\ast}=\emptyset, \ \forall
N_{i}\in\mathcal{R}.$$ Trivially it follows that $\bigcup_{n\in
J_{0}}N_{\overline{\omega},n}^{\ast}\subseteq \overline{N_{\omega}''}$ and
$\overline{N_{\omega}''}\cap U_{\omega}=\emptyset$. Let us recapitulate:
\begin{enumerate}   \item[1)] $N_{\omega}''\supseteq^{*} N_{\overline{\omega},n}$ for every $n\in
J_{0}$;\\
\item[2)] $|N_{\omega}''\cap N_{i}|<\omega$ for every
$i\in\omega$.
\end{enumerate}
Notice that $(N_{\omega}''')^{\ast}=(N_{\omega}')^{\ast}\cup
(N_{\omega}'')^{\ast}$ is a clopen subset of $\omega^{\ast}$. Now it turns out
that
\begin{displaymath}
(N_{\omega}''')^{\ast}\cap U_{\omega}= [(N_{\omega}')^{\ast}\cup
(N_{\omega}'')^{\ast}]\cap U_{\omega}= (j_{\alpha+1}^{\omega})^{-1}(\bigsqcup_{n\in
J_{1}} U_{n})\cup \emptyset=\bigsqcup_{n\in J_{1}} E_{n}^{\ast}
\end{displaymath}
and then we can conclude that $N_{\omega}'''\supseteq^{\ast}E_{n}$ for every $n\in J_{1}$
and $N_{\omega}'''\supseteq^{\ast}N_{\overline{\omega},n}$ for every $n\in J_{0}$. Let us
set
\begin{displaymath}
  M_{n}=\left\{
\begin{array}{ll}
     N_{\omega}'''\cap E_{n} \textmd{ if
} n\in J_{1} \\
      N_{\omega}'''\cap N_{\overline{\omega},n} \textmd{ if } n\in J_{0}.\\
\end{array}
\right. \end{displaymath}
 Certainly $M_{n}^{\ast}=E_{n}^{\ast}$
for every $n\in J_{1}$ and $M_{n}^{\ast}=N_{\overline{\omega},n}^{\ast}$ for every $n\in
J_{0}$.\\
For every $n\in\omega$ let us fix a point $l_{n}\in
M_{n}\backslash (\bigcup_{j=0}^{n-1}M_{j}\cup
\{l_{0},\ldots,l_{n-1}\})$ - it is possible since the family
$\{M_n\}$ is almost disjoint - and let us
set $L=\{l_{i}:i\in\omega\}$.\\
 Let us define
\begin{displaymath}
N_{\omega}=\bigsqcup_{n\in\omega} (M_{n}\backslash
\bigcup_{j=0}^{n-1}M_{j})\backslash\{l_{i}:i\in\omega\}=\bigsqcup_{n\in\omega}H_{n}
\end{displaymath}
where $H_{n}=(M_{n}\backslash \bigcup_{j=0}^{n-1}M_{j})\backslash\{l_{i}:i\in\omega\}$.
Notice that $N_{\omega}^{\ast}\supseteq \bigcup \mathcal{C}'$ (indeed from every $M_{n}$
      we removed only a finite number of points) and that
     $(N_{\omega}''')^{\ast}\backslash N_{\omega}^{\ast}\neq{\emptyset}$
      (since $N_{\omega}'''\backslash N_{\omega}= \{l_{i}:i\in\omega\}$)
       whence $|\omega\backslash N_{\omega}|=\omega$.\\
     Now let us take into account $N_{\omega}=\bigsqcup_{n\in\omega}H_{n}$
     where
$M_{n}^{\ast}=H_{n}^{\ast}$ for every $n\in\omega$; we want to remark that on each
$\overline{H_{n}}$ with $n\in J_{1}$ we have
already a decomposition of type $l(x_{n})$ by construction.\\
 Now
if $|J_{1}|=\omega$, on every $\overline{H_{n}}$ with $n\in J_{0}$
let us put a decomposition of type 1; then let us order the
subsets $\overline{H_{n}}$ in such a way that the types of
decomposition that we have put on them form a nondecreasing
sequence with the supremum equal to
$\beta_{\omega}<\alpha$.\\
 If
$|J_{1}|<\omega$, it turns out that $\sup\{l(x_{n}):n\in
J_{1}\}=\beta_{\omega}$ is a successor ordinal; let us put a
decomposition of type $\beta_{\omega}$ on every $\overline{H_{n}}$
with $n\in J_{0}$\footnote{In this case $|J_{0}|=\omega$.} and
then let us order the subsets $\overline{H_{n}}$ in such a way
that the types of decomposition that we have put on them form a
nondecreasing sequence
whose supremum is equal to $\beta_{\omega}<\alpha$.\\
 If
$J_{1}=\emptyset$, we choose to put a decomposition of type a
successor ordinal $\beta_{\omega}<\alpha$ on every
$\overline{H_{n}}$ with $n\in J_{0}$ and then we proceed as in the
latter case.\\ This time let us apply property S.6 to the subsets
$H_{n}$ and to $N_{\omega}$: it turns out that
$\{H_{n}\}_{n\in\omega}$ is a countably infinite family of
infinite pairwise disjoint subsets of $N_{\omega}$ and on every
$\overline{H_{n}}$ is given a relation of some type in such a way
that the supremum of the nondecreasing sequence consisting of the
types of decomposition has supremum $\beta_{\omega}$ with
$\beta_{\omega}$ that can take up value from $1$ to $\alpha$ not
included. Then the relation on
$\bigcup_{n=1}^{\infty}\overline{H_{n}}$
     obtained in this way can be extended to a relation $\approx_{\beta_{\omega}+1}$
     on $\overline{N_{\omega}}$ of type $\beta_{\omega}+1$ where $\beta_{\omega}+1$ can have value a successor
     ordinal from $2$ up to $\alpha$.\footnote{In the case in which $\alpha$
       is a limit ordinal, $\beta_{\omega}+1$ can take up value on the successor
     ordinals from $2$ up to an ordinal strictly less than $\alpha$.}

\begin{remark}\label{intvuota}We want to remark that for the
points constructed by the decompositions on the $\overline{H}_{n}$
with $n\in J_{0}$ it is always possible to find a fundamental
system of elementary neighborhoods contained in
$\overline{N}_{\omega}/\approx_{\beta_{\omega}+1}$ and such that
their inverse images through $j^{\omega}_{\alpha+1}$ have empty
intersection with $U_{\omega}$ since $H_{n}^{\ast}\cap
U_{\omega}=\emptyset$; from now on, we consider only these
neighborhoods as elementary neighborhoods of those points.
\end{remark}
Let us check the following properties:
\begin{enumerate}
    \item [T.1]\label{T.1} $\overline{N_{\omega}^{\ast}\backslash\bigcup_{i\in\omega}
    N_{i}^{\ast}}\neq \emptyset$:
    indeed it turns out that $(N_{\omega})^{\ast}
     \supseteq \bigcup \mathcal{L}_{0}$ while
      $\bigcup \mathcal{L}_{0}\cap
      (\bigcup_{i\in\omega}N_{i}^{\ast})=\emptyset$; hence, if $\mathcal{L}_{0}\neq\emptyset$,
       it follows that $N_{\omega}^{\ast}\backslash\bigcup
    N_{i}^{\ast}\supseteq \bigcup\mathcal{L}_{0}\neq\emptyset$.
    On the other hand if
     $\mathcal{L}_{0}=\emptyset$, then the family $\{H_{n}:n\in J_1\}$
      of pairwise disjoint subsets of $\omega$ is infinite; thus we can construct an
     infinite subset $T\subset N_{\omega}$ in this way: we choose a point $t_{m}\in
    H_{m}$ for every $m\in J_{1}$ and we set $T=\{t_{m}: m\in J_{1}\}$. The non-empty subset
     $T^{\ast}$ is such that $T^{\ast}\subseteq N_{\omega}^{\ast}$, while $T^{\ast}\cap
    N_{i}^{\ast}=\emptyset$ for every $i\in\omega$. We want to check it: certainly $T^{\ast}\cap
    H_{n}^{\ast}=\emptyset$ for every $n\in J_{1}$ and hence $T^{\ast}\cap
     (\bigsqcup H_{n}^{\ast})=\emptyset$; if there is an index $i\in\omega$ such that
      $|T\cap N_{i}|=\omega$, then $(T\cap N_{i})^{\ast}\subset (N_{\omega}^{\ast}\cap
      \bigcup_{i\in\omega}N_{i}^{\ast})$, while we know that
    $(N_{\omega}^{\ast}\cap \bigcup_{i\in\omega}N_{i}^{\ast})=\bigsqcup
    H_{n}^{\ast}$.

    \item [T.2] $\overline{\bigcup_{i\leq \omega}N^{\ast}_{i}}\neq \omega^{\ast}$: notice that the set $L$
    is such that
$L^{\ast}\cap
    N_{\omega}^{\ast}=\emptyset$ and that $L^{\ast}\cap
    N_{i}^{\ast}=\emptyset$ for every $i<\omega$ (indeed for
    every $i<\omega$ it holds that $|N_{i}\cap L|<\omega$; if there exists
    an index $i\in\omega$ such that $|L\cap N_{i}|=\omega$,
    then $(L\cap N_{i})^{\ast}\subset ((N_{\omega}''')^{\ast}\cap
    \bigcup_{i\in\omega}N_{i}^{\ast})$, while we know that
    $((N_{\omega}''')^{\ast}\cap \bigcup_{i\in\omega}N_{i}^{\ast})=\bigsqcup H_{n}^{\ast}$ and $L^{\ast}\cap (\bigsqcup
    H_{n}^{\ast})=\emptyset$).
    Then we obtain that $\bigcup_{i\leq\omega}N^{\ast}_{i}\cap
    L^{\ast}=\emptyset$ where
    $L^{\ast}$ is open in $\omega^{\ast}$ whence $\omega^{\ast}\backslash L^{\ast}$
     is a closed subset that contains
$\bigcup_{i\leq\omega}N^{\ast}_{i}$; so it contains its closure and it follows that
$\overline{\bigcup_{i\leq
    \omega}N^{\ast}_{i}}\cap L^{\ast}=\emptyset$. At the end, we can conclude that
$\overline{\bigcup_{i\leq
    \omega}N^{\ast}_{i}}\neq \omega^{\ast}$.
    \item [T.3] For every $i\in\omega$, the relations
    $\approx_{\beta_{i}+1}$ and $\approx_{\beta_{\omega}+1}$
    coincide on
    $\overline{N_{i}}\cap\overline{N_{\omega}}$: indeed $N_{\omega}^{\ast}
    \cap N_{i}^{\ast}\subseteq\bigsqcup H_{n}^{\ast}$ (with $n\in J_1$),
     the relation on $N_{\omega}^{\ast}$ extends the relations
      placed on the subsets $H_{n}^{\ast}$ (with $n\in J_1$) and these last relations
       coincide with the relations we put on the subsets
       $N_{i}^{\ast}$.\\ Then a relation $Q_{\omega+1}$ is defined on
$U_{\omega+1}=\bigcup_{i=1}^{\omega} N_{i}^{\ast}$.
 \item [T.4] A family
$\mathcal{C}_{\xi}\in\Gamma$ with index $\xi\leq \omega$ is not an
$(\omega+1)-$family. Notice that the families $\mathcal{C}_{\xi}$ with
$\xi<\omega$ do not exist and hence we have only to prove that the family $C_{\omega}$ is not an $(\omega+1)-$family.\\
\textit{We say that $\mathcal{C}_{\xi}\in{\Gamma}$ is an $(\omega+1)-$family if
$\mathcal{C}_{\xi}$ can be decomposed into two subfamilies
$\mathcal{L}_{0}^{\omega+1}$ and $\mathcal{L}_{1}^{\omega+1}$ satisfying the
following conditions.
\begin{enumerate}
\item [U.1] $\bigcup \mathcal{L}_{0}^{\omega+1}\cap
U_{\omega+1}=\emptyset$.
 \item [U.2] For every
$c\in{\mathcal{L}_{1}}^{\omega+1}$ there exists $i\leq\omega$, a point
$x_{c}\in{\overline{N_{i}}/\approx_{\beta_{i}+1}}$ of level
 $\gamma_{c}+1$ and an elementary neighborhood $U_{c}$ of $x_{c}$
 such that
 $c=j_{\beta_{i}+1}^{-1}(U_{c})\cap\omega^{\ast}$.
 \item [U.3] The set $\left\{x_{c}:
c\in{\mathcal{L}_{1}^{\omega+1}}\right\}$ is devoid of any accumulation point
in $U_{\omega+1}/Q_{\omega+1}$. \item [U.4] It holds that
$\sup\left\{\gamma_{c}+1:c\in \mathcal{L}_{1}^{\omega+1}\right\}<\alpha$.
\end{enumerate}}
Remember that $\mathcal{C}_{\overline{\omega}}$ is the $\omega$-family with
minimum index we have just used in order to construct
$\overline{N_{\omega}}/\approx_{\beta_{\omega}+1}$; moreover notice that if
$\overline{\omega}>\omega$ the family $\mathcal{C}_{\omega}$ is not an
$\omega$-family and then it neither is an $(\omega+1)$-family. Towards a
contradiction, suppose that it is an $(\omega+1)$-families: in
$\mathcal{L}^{\omega}_{0}$ we put
 the elements that lie in $\mathcal{L}^{\omega+1}_{0}$ and all
 those elements $c\in{\mathcal{L}^{\omega+1}_{1}}$ such that $\omega$ is the
 only value of the index
$i$ for which U.2 is satisfied; these $c$ are such that $c\cap U_{\omega}=\emptyset$ by
Remark \ref{intvuota}. Instead in $\mathcal{L}^{\omega}_{1}$ we put all the other
$c\in{\mathcal{L}^{\omega+1}_{1}}$ which are left: they obviously satisfy U.4 since we
are estimating the supremum on a lesser number of elements; moreover they satisfy U.3
since if the points $x_{c}$ had accumulation points in $U_{\omega}/Q_{\omega}$ then they
would have accumulation points in $U_{\omega+1}/Q_{\omega+1}$ due to the fact that the
new relation respects the old ones.
 \\
If instead $\overline{\omega}=\omega$, then $\mathcal{C}_{\overline{\omega}}=
\mathcal{L}_{0} \cup\mathcal L_{1}$ is not an $(\omega+1)$-family since the
elements of $\mathcal{C}_{\overline{\omega}}$ would have all to stay in
$\mathcal{L}_{1}^{\omega+1}$ but the corresponding infinite points $\{x_{c}:
c\in \mathcal{L}_{1}^{\omega+1}\}$, which are all in the compact space
$N_{\omega}^{\ast}/\approx_{\beta_{\omega}+1}$, must have an accumulation point
in $U_{\omega+1}/Q_{\omega+1}\supseteq
N_{\omega}^{\ast}/\approx_{\beta_{\omega}+1}$.

\end{enumerate}

\vspace{0,2cm}

 Suppose to have
constructed infinite subsets $N_{\gamma}\subset\omega$ (with
$\gamma<\delta<\omega_{1}$) and relations
$\approx_{\beta_{\gamma}+1}$ of type $\beta_{\gamma}+1<\alpha+1$
on $\overline{N_{\gamma}}$ satisfying the following conditions.
\begin{enumerate}
    \item [T.1] $\overline{N_{\gamma}^{\ast}\backslash\bigcup_{\gamma'<\gamma}N_{\gamma'}^{\ast}}\neq
    \emptyset$.
    \item [T.2] $\overline{\bigcup_{\gamma'\leq \gamma}N^{\ast}_{\gamma'}}\neq \omega^{\ast}$.
    \item [T.3] For every $\gamma'<\gamma$, the relations
    $\approx_{\beta_{\gamma'}+1}$ and $\approx_{\beta_{\gamma}+1}$
    coincide on $\overline{N_{\gamma'}}\cap\overline{N_{\gamma}}$.
    \item [T.4] A family $\mathcal{C}_{\xi}$ with index
     $\xi\leq \gamma$ is not a $(\gamma+1)-$family.
\end{enumerate}
By T.3 and the decompositions set on the $\overline{N}_{\gamma}$
for every $\gamma<\delta$, a decomposition $Q_{\delta}$ is defined
on $U_{\delta}=\bigcup_{\gamma<\delta}N_{\gamma}^{\ast}$. Let
$j^{\delta}_{\alpha+1}$ be the quotient mapping
$j^{\delta}_{\alpha+1}:U_{\delta}\rightarrow
U_{\delta}/Q_{\delta}$.\\ We say that an element
$\mathcal{C}_{\xi}\in{\Gamma}$ is a $\delta-$family if
$\mathcal{C}_{\xi}$ can be decomposed into two subfamilies
$\mathcal{L}_{0}$ and $\mathcal{L}_{1}$ satisfying the following
conditions.
\begin{enumerate}
\item [U.1] $\bigcup \mathcal{L}_{0}\cap U_{\delta}=\emptyset$.
 \item [U.2] For every
$c\in{\mathcal{L}_{1}}$ there exists $\gamma<\delta$, a point
$x_{c}\in{\overline{N_{\gamma}}/\approx_{\beta_{\gamma}+1}}$ of
level
 $\gamma_{c}+1$ and an elementary neighborhood $U_{c}$ of $x_{c}$
 such that
 $c=j_{\beta_{\gamma}+1}^{-1}(U_{c})\cap\omega^{\ast}$.
 \item[U.3]
The set $\left\{x_{c}: c\in{\mathcal{L}_{1}}\right\}$ has no
accumulation point in $U_{\delta}/Q_{\delta}$. \item [U.4] It
holds that
$\sup\left\{\gamma_{c}+1:c\in{\mathcal{L}_{1}}\right\}<\alpha$.
\end{enumerate}
We can rewrite these properties in the following way:
\begin{enumerate} \item[i)] $\mathcal{L}_{0}$ has to consist
 of elements $C_{n}^{\ast}$ where the subsets $C_{n}$
are transversal to the subsets $N_{\gamma}$, i.e. every
$C_{n}\subset\omega$ intersects every $N_{\gamma}$ in a finite
number of points (in this way we are respecting $U.1$); \item[ii)]
$\mathcal{L}_{1}$ has to consist of elements $C_{m}^{\ast}$ where
for every $m$ there exists $\gamma\in\delta$, a point $x_{m}\in
\overline{N_{\gamma}}/\approx_{\beta_{\gamma}+1}$ of level
$l(x_{m})<\alpha$ and an elementary neighborhood $U(x_{m})$ such
that
$C_{m}^{\ast}=j_{\beta_{\gamma}+1}^{-1}(U(x_{m}))\cap\omega^{\ast}$.
A further necessary request is that the set $\{x_{m}\}$ is devoid
of any accumulation point in $U_{\delta}/Q_{\delta}$ and that
$\sup\{l(x_{m})\}<\alpha$. We want to remark that by $l(x_{m})$ we
mean a successor ordinal. (In this way we are respecting
U.2-U.3-U.4)
 \end{enumerate}
Let us show that it is possible to find a $\delta$-family
$C_{\xi}$: for example, we can use again Lemma \ref{lemma01},
since the subsets $N_{\gamma}^{\ast}$ with
$\gamma<\delta<\omega_{1}$ comply with the hypotheses. We want to
remark that here the fact that $\delta$ is a countable ordinal is
essential in order to apply Lemma \ref{lemma01}. In this way we
are able to find a subset $\Delta_{\delta}\subset\omega$ with
$|\Delta_{\delta}|=\omega$ and such that $\Delta_{\delta}$
intersects every $N_{\gamma}$ in a finite number of points. We can
decompose this infinite set in an infinite number of infinite
subsets $T_{n}\subset\omega$ which again intersect every
$N_{\gamma}$ in a finite number of points; we set
$\mathcal{L}_{0}=\{T_n^{\ast}:n\in\omega\}$. It is clear that
$\mathcal{L}=\mathcal{L}_{0}$ is a $\delta$-family.\\ Among all
the $\delta$-families in $\Gamma$, let us take the one with the
minimum index $\overline{\delta}$: we can write it as
$\mathcal{C}_{\overline{\delta}}=\mathcal{L}_{0}\sqcup\mathcal{L}_{1}$
where $\mathcal{L}_{0}=\{N_{\overline{\delta},n}^{\ast}:n\in
J_0\}$, $\mathcal{L}_{1}=\{N_{\overline{\delta},n}^{\ast}:n\in
J_1\}$ and $J_0\cap J_1=\emptyset$. Of course, by construction,
the $\delta$-family $\mathcal{C}_{\overline{\delta}}$ will comply
with the following properties.
\begin{enumerate}
\item [U.1] $\bigcup \mathcal{L}_{0}\cap U_{\delta}=\emptyset$;
 \item [U.2] For every $N_{\overline{\delta},n}^{\ast}$
 with $n\in J_{1}$ there exist an index $\gamma_{n}\in\delta$,
 a point $x_{n}
  \in \overline{N_{\gamma_{n}}}/\approx_{\beta_{\gamma_{n}}+1}$
   of level $l(x_{n})<\alpha$
   and an elementary neighborhood
   $U(x_{n})$ such that
   $$N_{\overline{\delta},n}^{\ast}=j_{\beta_{\gamma_{n}}+1}^{-1}(U(x_{n}))
   \cap\omega^{\ast}.$$
 \item [U.3] The set $\{x_{n}:n\in
J_{1}\}$ has no accumulation point in $U_{\delta}/Q_{\delta}$.
\item [U.4] It holds that $\sup\{l(x_{n}):n\in
J_{1}\}=\beta_{\delta}<\alpha$ with $\beta_{\delta}$ that can take
up value from $1$ to $\alpha$ not included. Without loss of
generality, we can always assume that the levels of the points are
ordered in a nondecreasing way.
\end{enumerate}
For every $n\in J_{1}$ it holds that $\hat{U}(x_{n})=U(x_{n})\backslash\omega$
is a clopen neighborhood of $x_{n}$ in
$N_{\gamma_{n}}^{\ast}/\approx_{\beta_{\gamma_{n}}+1}$. In order to apply Lemma
\ref{lemma1}, we need to rewrite $\bigcup_{\gamma\in\delta} N_{\gamma}^{\ast}$
as a disjoint union; notice that, since $\delta$ is a countable ordinal, we can
enumerate $\{N_{\gamma}^{\ast}\}_{\gamma\in\delta}$ as
$\{N_{\gamma_i}^{\ast}\}_{i\in\omega}$. Let us set
$\tilde{N}_{\gamma_1}=N_{\gamma_1}, \tilde{N}_{\gamma_2}=N_{\gamma_2}\backslash
N_{\gamma_1},\ldots, \tilde{N}_{\gamma_k}=N_{\gamma_k}\backslash
\bigcup_{i=1}^{k-1}N_{\gamma_{i}}$. It holds that
$\bigsqcup_{k\in\omega}\tilde{N}_{\gamma_k}=\bigcup_{k\in\omega}N_{\gamma_k}$,
whence $\bigsqcup_{k\in\omega}\tilde{N}_{\gamma_k}^{\ast}=
\bigcup_{k\in\omega}N_{\gamma_k}^{\ast}$; we have only to prove the non-trivial
inclusion $\bigsqcup_{k\in\omega}\tilde{N}_{\gamma_k}^{\ast}\supseteq
\bigcup_{k\in\omega}N_{\gamma_k}^{\ast}$: if
$x\in\bigcup_{k\in\omega}N_{\gamma_k}^{\ast}$, then there is
$\bar\imath\in\omega$ such that $$x\in
N_{\gamma_{\bar\imath}}^{\ast}=\big[\bigsqcup_{k\leq\bar\imath}\tilde{N}_{\gamma_k}\big]^{\ast}=
\bigsqcup_{k\leq\bar\imath}\tilde{N}_{\gamma_k}^{\ast}\subseteq\bigsqcup_{k\in\omega}\tilde{N}_{\gamma_k}^{\ast}.$$
Notice that, for every $k\in\omega$, $\tilde{N}_{\gamma_k}^{\ast}$ is
distinguished relative to $Q_{\delta}$. Finally we can apply Lemma \ref{lemma1}
since the countable family $\{\tilde{N}_{\gamma}\}_{\gamma<\delta}$, the points
$\{x_{n}\}_{n\in J_{1}}$ and the relation $Q_{\delta}$ defined on
$U_{\delta}=\bigsqcup_{\gamma<\delta}
\tilde{N}_{\gamma}^{\ast}=\bigcup_{\gamma<\delta}N_{\gamma}^{\ast}$ satisfy the
hypotheses.\footnote{At most we have to restrict the neighborhoods of the
points $\{x_n\}$ in such a way that each of them belongs to some
$\tilde{N}_{\gamma}^{\ast}/Q_{\delta}$ for some $\gamma<\delta$.} We remark
that $\bigcup\hat{U}(x_{n})\neq U_{\delta}/Q_{\delta}$ since in
$U_{\delta}/Q_{\delta}$ there are points of level $\alpha$
 that $\bigcup\hat{U}(x_{n})$ does not cover.\footnote{Notice that $\bigcup\hat{U}(x_{n})\neq
U_{\delta}/Q_{\delta}$ also in the case in which $\alpha$ is a limit ordinal:
 indeed at the beginning of the construction we put decompositions of type $\beta_{i}+1$ on
 the subsets $\overline{N_{i}}$ in such a way that
 $\sup\{\beta_{i}+1\}=\alpha$; hence in $U_{\delta}/Q_{\delta}$ there
 certainly exists a point of level $\beta_{\delta}+1<\alpha$ that $\bigcup\hat{U}(x_{n})$ does not cover.}\\
Therefore it is possible to find pairwise disjoint elementary neighborhoods
$U_{n}$ with $x_{n}\in U_{n}$ and a subset $N_{\delta}'\subset \omega$ such
that $$(N_{\delta}')^{\ast}\cap U_{\delta}=\bigsqcup
(j_{\alpha+1}^{\delta})^{-1}(U_{n})=\bigsqcup E_{n}^{\ast}.$$
 Let us define
$\mathcal{C}'=\mathcal{L}_{0}\cup
\{(j_{\alpha+1}^{\delta})^{-1}(U_{n}):n\in J_{1}\}$.\\
 Let us set
$\mathcal{Q}=\{N_{\overline{\delta},n}: n\in J_{0}\}$ and $\mathcal{R}=\{N_{\gamma}:
\gamma\in\delta\}$. Then $\mathcal{P}=\mathcal{Q}\cup\mathcal{R}$ is a family with the
following properties:
\begin{itemize}
 \item[-] $\mathcal{Q}$ is an almost disjoint family;
\item[-] $|\mathcal{Q}|\leq\omega$ and $|\mathcal{R}|\leq\omega$;
\item[-] for every $N_{\overline{\delta},n}\in\mathcal{Q}$ and
every $N_{\gamma}\in\mathcal{R}$ it holds that $|N_{\overline{\delta},n}\cap
N_{\gamma}|<\omega.$
\end{itemize}
  Therefore by Lemma \ref{lemma2}
  there exists a subset $N_{\delta}''\in[\omega]^{\omega}$
   such that $$\bigcup_{n\in J_{0}}N_{\overline{\delta},n}^{\ast}\subseteq
(N_{\delta}'')^{\ast}\quad \textmd{and}\quad (N_{\delta}'')^{\ast}\cap
 N_{\gamma}^{\ast}=\emptyset,\quad \forall N_{\gamma}\in\mathcal{R};$$
 obviously it holds that $\bigcup_{n\in J_{0}}N_{\overline{\delta},n}^{\ast}\subseteq \overline{N_{\delta}''}$ and
$\overline{N_{\delta}''}\cap U_{\delta}=\emptyset$.  Hence
$N_{\delta}''$ is such that:
\begin{enumerate}   \item[1)] $N_{\delta}''\supseteq^{*}
 N_{\overline{\delta},n}$ for every $n\in
J_{0}$;\\
\item[2)] $|N_{\delta}''\cap N_{\gamma}|<\omega$ for every
$\gamma\in\delta.$
\end{enumerate}

Let us remark that $(N_{\delta}''')^{\ast}=(N_{\delta}')^{\ast}\cup
(N_{\delta}'')^{\ast}$ is a clopen subset of $\omega^{\ast}$. Now it turns out
that
\begin{displaymath}
(N_{\delta}''')^{\ast}\cap U_{\delta}= [(N_{\delta}')^{\ast}\cup
(N_{\delta}'')^{\ast}]\cap U_{\delta}= [(j_{\alpha+1}^{\delta})^{-1}(\bigsqcup_{n\in
J_{1}} U_{n})]\cup\emptyset =\bigsqcup_{n\in J_{1}} E_{n}^{\ast}
\end{displaymath}
and hence $N_{\delta}'''\supseteq^{\ast}E_{n}$ for every $n\in J_{1}$ and
$N_{\delta}'''\supseteq^{\ast}N_{\overline{\delta},n}$ for every $n\in J_{0}$. Let us set
 $$\qquad M_{n}=\left\{
\begin{array}{ll}
     N_{\delta}'''\cap E_{n} \textmd{ if
} n\in J_{1} \\
      N_{\delta}'''\cap N_{\overline{\delta},n} \textmd{ if } n\in J_{0}.\\
\end{array}
\right.$$
 Certainly $M_{n}^{\ast}=E_{n}^{\ast}$ for every $n\in
J_{1}$ and $M_{n}^{\ast}=N_{\overline{\delta},n}^{\ast}$ for every $n\in J_{0}$.\\ For
every $n\in\omega$ let us fix a point $l_{n}\in M_{n}\backslash
(\bigcup_{j=0}^{n-1}M_{j}\cup \{l_{0},\ldots,l_{n-1}\})$ and let us set
$L=\{l_{i}:i\in\omega\}$. Let us define
\begin{displaymath}
N_{\delta}=\bigsqcup_{n\in\omega} (M_{n}\backslash
\bigcup_{j=0}^{n-1}M_{j})\backslash\{l_{i}:i\in\omega\}=\bigsqcup_{n\in\omega}H_{n}.
\end{displaymath}
Notice that $N_{\delta}$ is such that
     $N_{\delta}^{\ast}\supseteq \bigcup \mathcal{C}'$  (indeed from every $M_{n}$
      we removed only a finite number of points) and at the same
      time that
     $(N_{\delta}''')^{\ast}\backslash N_{\delta}^{\ast}\neq{\emptyset}$
      (since $N_{\delta}'''\backslash N_{\delta}= \{l_{i}:i\in\omega\}$)
       whence $|\omega\backslash N_{\delta}|=\omega$.\\
      Therefore it holds that $N_{\delta}=\bigsqcup_{n\in\omega}H_{n}$, with
$M_{n}^{\ast}=H_{n}^{\ast}$ for every $n\in\omega$. Remember that
on each $\overline{H_{n}}$ with $n\in J_{1}$ we have already a
decomposition of type $l(x_{n})$.\\ If $|J_{1}|=\omega$, on every
$\overline{H_{n}}$ with $n\in J_{0}$ let us set a decomposition of
type 1 and let us order the subsets $\overline{H_{n}}$ in such a
way that the types of decomposition that we have put on them form
a nondecreasing sequence with supremum equal to
$\beta_{\delta}<\alpha$.\\ If $|J_{1}|<\omega$, it turns out that
$\sup\{l(x_{n}):n\in J_{1}\}=\beta_{\delta}$ is a successor
ordinal; let us set decompositions of type $\beta_{\delta}$ on
every $\overline{H_{n}}$ with $n\in J_{0}$\footnote{In this case
$|J_{0}|=\omega$} and let us order the subsets $\overline{H_{n}}$
in such a way that the types of decomposition that we have put on
them form a nondecreasing sequence with supremum equal to
$\beta_{\delta}<\alpha$.\\
 If
$J_{1}=\emptyset$, we choose to put a decomposition of type a successor ordinal
$\beta_{\delta}<\alpha$ on every $\overline{H_{n}}$ with $n\in J_{0}$ and then we proceed
as in the
latter case.\\
 Let us apply property S.6 to the subsets
$H_{n}$ and to $N_{\delta}$: indeed $\{H_{n}\}_{n\in\omega}$ is a
countably infinite family of infinite pairwise disjoint subsets of
$N_{\delta}$ and on every $\overline{H_{n}}$ is given a relation
of some type in such a way that the supremum of the nondecreasing
sequence consisting of the types of decomposition has supremum
$\beta_{\delta}$ with $\beta_{\delta}$ that can take up value from
$1$ to $\alpha$ not included. Then the relation on
$\bigcup_{n=1}^{\infty}\overline{H_{n}}$
     obtained in this way can be extended to a relation $\approx_{\beta_{\delta}+1}$
     on $\overline{N_{\delta}}$ of type $\beta_{\delta}+1$ with $2\leq\beta_{\delta}+1\leq\alpha$.\footnote{In the case in which $\alpha$
       is a limit ordinal it turns out that $2\leq\beta_{\delta}+1<\alpha$}\\
\begin{remark}\label{intvuotadelta}We want to remark that for the
points constructed by the decompositions on the $\overline{H}_{n}$
with $n\in J_{0}$ it is always possible to find a fundamental
system of elementary neighborhoods contained in
$\overline{N}_{\delta}/\approx_{\beta_{\delta}+1}$ and such that
their inverse images through $j^{\delta}_{\alpha+1}$ have empty
intersection with $U_{\delta}$ since $H_{n}^{\ast}\cap
U_{\delta}=\emptyset$; from now on, we consider only these
neighborhoods as elementary neighborhoods of those points.
\end{remark}

Let us check that properties T.1 to T.4 hold.
\begin{enumerate}
    \item [T.1] $\overline{N_{\delta}^{\ast}\backslash(\bigcup_{\gamma\in\delta}
    N_{\gamma}^{\ast})}=\overline{N_{\delta}^{\ast}\backslash(\bigcup_{\gamma\in\delta}
    \tilde{N}_{\gamma}^{\ast})}\neq \emptyset$:
    the check is the same as in the case of
    $N_{\omega}$ (see page \pageref{T.1}).

    \item [T.2] $\overline{\bigcup_{\gamma\leq \delta}N^{\ast}_{\gamma}}\neq \omega^{\ast}$:
   the check of this property is similar to that we have just done in the
   case of $N_{\omega}$ (see page \pageref{T.1}).
    \item [T.3] For every $\gamma\in\delta$, the relations
    $\approx_{\beta_{\gamma}+1}$ and $\approx_{\beta_{\delta}+1}$ coincide on
    $\overline{N_{\gamma}}\cap\overline{N_{\delta}}$: indeed
     $N_{\gamma}^{\ast}\cap N_{\delta}^{\ast}\subseteq U_{\delta}\cap N_{\delta}^{\ast}\subseteq
    \bigcup H_{n}^{\ast}$ (with $n\in J_1$),
     the relation on $N_{\delta}^{\ast}$ extends the relations
      set on the subsets $H_{n}^{\ast}$
      and these last relations
       coincide with those we put on the subsets $N_{\gamma}^{\ast}$ by construction.
Then a relation $Q_{\delta+1}$ is defined on
$U_{\delta+1}=\bigcup_{\gamma=1}^{\delta} N_{\gamma}^{\ast}$.
\item [T.4] A family $\mathcal{C}_{\xi}$ with index $\xi\leq
\delta$ is not a $(\delta+1)-$family.\\ \textit{We say that
$\mathcal{C}_{\xi}\in{\Gamma}$ is a $(\delta+1)-$family if
$\mathcal{C}_{\xi}$ can be decomposed into two subfamilies
$\mathcal{L}_{0}^{\delta+1}$ and $\mathcal{L}_{1}^{\delta+1}$
satisfying the following conditions.
\begin{enumerate}
\item [U.1] $\bigcup \mathcal{L}_{0}^{\delta+1}\cap
U_{\delta+1}=\emptyset$;
 \item [U.2] For every
$c\in{\mathcal{L}_{1}}^{\delta+1}$ there exist $\gamma\leq\delta,$ a point
$x_{c}\in{\overline{N_{\gamma}}/\approx_{\beta_{\gamma}+1}}$ of level
 $\gamma_{c}+1$ and an elementary neighborhood $U_{c}$ of $x_{c}$
 such that
 $c=j_{\beta_{\gamma}+1}^{-1}(U_{c})\cap\omega^{\ast}$.
 \item [U.3] The set $\left\{x_{c}:
c\in{\mathcal{L}_{1}^{\delta+1}}\right\}$ is devoid of any
accumulation point in $U_{\delta+1}/Q_{\delta+1}$. \item [U.4] It
holds that $\sup\left\{\gamma_{c}+1:c\in
\mathcal{L}_{1}^{\delta+1}\right\}<\alpha$.
\end{enumerate}}
Remember that $\mathcal{C}_{\overline{\delta}}$ is the $\delta$-family with minimum index
we have just used in the construction of
$\overline{N}_{\delta}/\approx_{\beta_{\delta}+1}$. If $\overline{\delta}>\delta$ the
families $\mathcal{C}_{\xi}$ with $\xi\leq\delta$ are not $\delta$-families and then they
neither are $(\delta+1)$-families. Towards a contradiction, suppose that they are
$(\delta+1)$-families; then in $\mathcal{L}^{\delta}_{0}$ we put
 the elements that lie in $\mathcal{L}^{\delta+1}_{0}$ and all
 those elements $c\in{\mathcal{L}^{\delta+1}_{1}}$ such that $\delta$ is the
 only value of the index
$\gamma$ for which U.2 is satisfied; by Remark \ref{intvuotadelta}
these $c$ are such that $c\cap U_{\delta}=\emptyset$. Instead in
$\mathcal{L}^{\delta}_{1}$ we put all the other
$c\in{\mathcal{L}^{\delta+1}_{1}}$ that are left: they obviously
satisfy U.4, since we are estimating the supremum on a lesser
number of elements; moreover they satisfy U.3, since if the points
$x_{c}$ had accumulation points in $U_{\delta}/Q_{\delta}$, then
they would have accumulation points in
$U_{\delta+1}/Q_{\delta+1}$, due to the fact that the new relation
respects the old ones.
 \\  On the other hand, if
$\overline{\delta}=\delta$, then the families $\mathcal{C}_{\xi}$ with $\xi<\delta$ are
not $\delta$-families and by what we have just remarked they neither are
$(\delta+1)$-families; on the other hand $\mathcal{C}_{\overline{\delta}}=
\mathcal{L}_{0}^{\delta} \cup\mathcal L_{1}^{\delta}$ is not a $(\delta+1)$-family, since
the elements of ${C}_{\overline{\delta}}$ would have all to stay in
$\mathcal{L}_{1}^{\delta+1}$ but the corresponding infinite points $x_{c}$, which are all
in the compact space $N_{\delta}^{\ast}/\approx_{\beta_{\delta}+1}$, must have an
accumulation point in $U_{\delta+1}/Q_{\delta+1}\supseteq
N_{\delta}^{\ast}/\approx_{\beta_{\delta}+1}$.
\end{enumerate}


Therefore, by transfinite induction, we have defined a relation $Q_{\omega_{1}}$ on
$\bigcup_{\gamma<\omega_{1}}\overline{N_{\gamma}}$ which coincides
with $\approx_{\beta_{\gamma}+1}$ on each $\overline{N_{\gamma}}$.\\
Let us prove the following lemma.

\begin{lemma}\label{osimp}
If a family $\mathcal{C}_{\xi}\in \Gamma$ is not a $\vartheta$-family then it is not a
$\delta$-family for every $\delta>\vartheta$.
\end{lemma}
\iniziodim We prove that if $\mathcal{C}_{\xi}\in \Gamma$
 is a $\delta$-family then it is also a $\vartheta$-family.
 Let us suppose that $\mathcal{C}_{\xi}$
is a $\delta$-family; then it can be decomposed into two
subfamilies $\mathcal{L}^{\delta}_{0}$ and
$\mathcal{L}^{\delta}_{1}$ satisfying the following conditions.
\begin{enumerate}
\item [U.1] $\bigcup \mathcal{L}^{\delta}_{0}\cap
U_{\delta}=\emptyset$.
 \item [U.2] For every
$c\in{\mathcal{L}^{\delta}_{1}}$ there exist $\gamma<\delta$, a point $x_{c}\in
\overline{N_{\gamma}}/\approx_{\beta_{\gamma}+1}$ of level
 $\gamma_{c}+1$ and an elementary neighborhood $U_{c}$ of $x_{c}$
 such that
  $c=j_{\beta_{\gamma}+1}^{-1}(U_{c})\cap\omega^{\ast}$.
\item [U.3] The set $\left\{x_{c}:
c\in{\mathcal{L}^{\delta}_{1}}\right\}$ has no accumulation point
in $U_{\delta}/Q_{\delta}$. \item [U.4] It holds that
$\sup\left\{\gamma_{c}+1:c\in{\mathcal{L}^{\delta}_{1}}\right\}<\alpha$.
\end{enumerate}

We want to show that $\mathcal{C}_{\xi}$ is also a $\vartheta$-family. In
$\mathcal{L}^{\vartheta}_{0}$ we put
 the elements that lie in $\mathcal{L}^{\delta}_{0}$ and all
 those elements $c\in{\mathcal{L}^{\delta}_{1}}$ for which the only
 ordinals that fit for U.2 are larger than or equal to $\vartheta$;
   these $c$ are the inverse images of elementary neighborhoods of points constructed by starting
    from some $\overline{F}_{n}$ where $F_{n}^{\ast}\in \mathcal{L}_{0}^{\zeta}$
     with $\zeta\geq\vartheta$.
      We know that the elementary neighborhoods of these points
       are contained in $\overline{F}_{n}/\approx_{\beta_{\zeta}+1}$ and then
        by Remarks \ref{intvuota} and \ref{intvuotadelta} it
        follows that for each of these $c$ it holds that
         $c\cap U_{\vartheta}=\emptyset$. On the other hand in $\mathcal{L}^{\vartheta}_{1}$
we put all the other $c\in{\mathcal{L}^{\delta}_{1}}$ that are left. They
obviously satisfy U.4, since we are estimating the supremum on a lesser number
of elements; moreover they satisfy U.3, since if the points $x_{c}$ had
accumulation points in $U_{\vartheta}/Q_{\vartheta}$, then they would have
accumulation points in $U_{\delta}/Q_{\delta}$ due to the fact that the new
relations respect the old ones. \finedim

Now we can state the following fundamental remark.

\begin{remark}\label{remark} For every $\vartheta<\omega_{1}$,
 $\mathcal{C}_{\vartheta}\in\Gamma$ is not an $\omega_{1}$-family.
 Towards a contradiction suppose that there exists an index $\vartheta<\omega_{1}$
 such that
$\mathcal{C}_{\vartheta}$ is an
 $\omega_{1}$-family. By transfinite induction we proved that,
for every $\vartheta<\omega_{1}$, a family
 $\mathcal{C}_{\xi}$ with $\xi\leq\vartheta$ is not a
 $(\vartheta+1)$-family and hence $\mathcal{C}_{\vartheta}$ is not a
  $(\vartheta+1)$-family. On the other hand we supposed that $\mathcal{C}_{\vartheta}$
  is an $\omega_{1}$-family and then by Lemma \ref{osimp}
   it is a $(\vartheta+1)$-family. A contradiction. Let us point out that there
    can not exist $\omega_{1}$-families in $\Gamma$, since the
elements of the set $\Gamma$
   have indeces that go from $\omega$ included to $\omega_{1}$ not included.
\end{remark}

Finally we define the relation $\approx_{\alpha+1}$ on
$\beta\omega$ in this way:
\begin{enumerate} \item[-] it coincides with $Q_{\omega_{1}}$ on
$\bigcup_{\gamma<\omega_{1}}\overline{N_{\gamma}}$; \item[-] two free ultrafilter
belonging to
   $\omega^{\ast}\backslash\bigcup_{\gamma<\omega_{1}}N_{\gamma}^{\ast}$
     are equivalent under the relation $\approx_{\alpha+1}$.
     \end{enumerate}
Let us call $K_{\alpha+1}$ the space obtained by the quotient of $\beta\omega$ with this
relation and $j_{\alpha+1}$ the natural quotient mapping. Let us remark that, by property
U.2, $\omega^{\ast}\backslash\bigcup_{\gamma<\omega_{1}}N_{\gamma}^{\ast}$ is not empty:
indeed for every
  $\gamma\in\omega_{1}$ it holds that
  $B_{\gamma}=\omega^{\ast}\backslash\bigcup_{\gamma'\leq\gamma}N_{\gamma'}^{\ast}$
  is a closed subset of $\omega^{\ast}$ and the subsets $B_{\gamma}$ (with $\gamma\in\omega_1$) are such
  that $B_{\gamma_{1}}\supseteq B_{\gamma_{2}}$ for every
  $\gamma_{1}<\gamma_{2}$; moreover the family of closed subsets $\{B_{\gamma}\}_{\gamma\in\omega_{1}}$
   has the finite intersection property by T.2 proved
    for every step $\gamma\in\omega_{1}$. Thus, due to the compactness
    of $\omega^{\ast}$, it follows that
   $$\bigcap_{\gamma\in\omega_{1}}B_{\gamma}=
   \bigcap_{\gamma\in\omega_{1}}(\omega^{\ast}\backslash
   \bigcup_{\gamma'\leq\gamma}N_{\gamma'}^{\ast})=\omega^{\ast}
   \backslash\bigcup_{\gamma\in\omega_{1}}N_{\gamma}^{\ast}\neq\emptyset.$$\\
Then $j_{\alpha+1}$ collapses
$\omega^{\ast}\backslash\bigcup_{\gamma<\omega_{1}}N_{\gamma}^{\ast}$
to a single point which we call $x_{\infty}$.\\ If an element of
$K_{\alpha+1}$ is a point of the decompositions
$\approx_{\beta_{\delta}+1}$ and $\approx_{\beta_{\gamma}+1}$,
then the point lies in the same level in
$\overline{N_{\delta}}/\approx_{\beta_{\delta}+1}$ and
$\overline{N_{\gamma}}/\approx_{\beta_{\gamma}+1}$: indeed
whenever we reconsider a point that was in some previous
decompositions we take care that there exists an elementary
neighborhood of it that accompanies the point in the new
decomposition; in this way the level of the point is preserved and
the definition of $L_{\beta+1}$ as the set of the points that lie
in the level $\beta+1$ in some
$\overline{N_{\gamma}}/\approx_{\beta_{\gamma}+1}$ is correct. If
a point of the space $K_{\alpha+1}$ is a point of the
decompositions $\approx_{\beta_{\delta}+1}$ and
$\approx_{\beta_{\gamma}+1}$ with $\delta>\gamma$, then the
problem reduces to examine what happens in
$\overline{N_{\delta}}/\approx_{\beta_{\delta}+1}$ as regards its
elementary neighborhoods.
 We have to remark that in the
construction of the space $K_{\alpha+1}$ we paid attention to the fact that for every
level $0<\beta+1\leq\alpha$\footnote{strictly smaller than $\alpha$ in the case in which
$\alpha$ is a limit ordinal.} every point of level $\beta+1$ had a basis of clopen
subsets homeomorphic
 to the space $K_{\beta+1}$ which is compact and sequential by inductive hypothesis.\\
Now we have to understand which are the elementary neighborhoods
of the unique point of level $\alpha+1$ in
  $K_{\alpha+1}$, i.e. of the point $x_{\infty}=j_{\alpha+1}(\omega^{\ast}\backslash\bigcup_{\gamma<\omega_{1}}N_{\gamma}^{\ast})$.
On this subject let us prove the following lemma.
 \begin{lemma}\label{lemma intorni punti con alpha successore}
  The collection of the clopen subsets $K_{\alpha+1}\backslash \bigcup_{x\in G}U_{x}$
   (where $G$ is a finite set and for every $x\in G$ the clopen subset $U_{x}$
 is an elementary neighborhood in $K_{\alpha+1}$ of the point $x$
  that can have level equal to a successor ordinal smaller than or equal
 to $\alpha$)
  is a basis at the point $x_{\infty}$.
\end{lemma}

\iniziodim In an obvious way $K_{\alpha+1}\backslash \bigcup_{x\in G}U_{x}$
 is a clopen subset of $K_{\alpha+1}$ containing $x_{\infty}$.
Let $A$ be an open subset of $K_{\alpha+1}$ containing $x_{\infty}$ and let
$C=K_{\alpha+1}\backslash A$ be the complementary closed subset. For every $x\in C$, let
$U_{x}$ be an elementary clopen neighborhood of $x$; trivially, by taking all the clopen
neighborhoods $U_{x}$, with $x\in C$, we cover $C$. Let us consider
$j_{\alpha+1}^{-1}(C)$: it is a closed subset of $\beta\omega$ and then it is compact. If
we take all the open subsets $j_{\alpha+1}^{-1}(U_{x})$ (with $x\in C$) they form an open
cover of $j_{\alpha+1}^{-1}(C)$; then there exists a finite subcover $\bigcup_{x\in
G}j_{\alpha+1}^{-1}(U_{x})\supseteq j_{\alpha+1}^{-1}(C)$. Hence it turns out that
$$j_{\alpha+1}(\bigcup_{x\in
G}j_{\alpha+1}^{-1}(U_{x}))=\bigcup_{x\in
G}j_{\alpha+1}(j_{\alpha+1}^{-1}(U_{x}))=$$$$\bigcup_{x\in
G}U_{x}\supseteq j_{\alpha+1}(j_{\alpha+1}^{-1}(C))=C$$ and, by
passing to the complementary subsets, we can conclude that
$K_{\alpha+1}\backslash \bigcup_{x\in G}U_{x} \subseteq
K_{\alpha+1}\backslash C=A$.
    \finedim
We call \emph{elementary} each of these neighborhoods of the point $x_{\infty}$.

\section{Check of the properties of $K_{\alpha+1}$}

Now we want to check that the space $K_{\alpha+1}$ satisfies all
the  requested properties.

\begin{lemma}
$K_{\alpha+1}$ is a Hausdorff space and it is compact.
\end{lemma}
\iniziodim Trivially the points of $L_{0}$ can be separated from every other
point since they are isolated. Moreover, if we want to separate
$x_{\infty}=j_{\alpha+1}(\omega^{\ast}
  \backslash\bigcup_{\gamma<\omega_{1}}
  N_{\gamma}^{\ast})$ from
any other point $x$, it is enough to take respectively the open
disjoint elementary neighborhoods $K_{\alpha+1}\backslash U_{x}$
and
    $U_{x}$.\\
Suppose now to have to part two points $x_{1}$ and $x_{2}$ of
level smaller than $\alpha+1$; it is possible to face up with two
different situations.
\begin{enumerate}
\item[1)] There exists $\vartheta\in\omega_{1}$ such that
$x_{1},x_{2}\in j_{\alpha+1}(\overline{N_{\vartheta}})$; notice
that $j_{\alpha+1}(\overline{N_{\vartheta}})\simeq K_{\beta+1}$
(with $\beta+1<\alpha+1$) which is a Hausdorff space by inductive
hypothesis. Then in $j_{\alpha+1}(\overline{N_{\vartheta}})$ there
are two open neighborhoods $V_{x_{1}}$ and $V_{x_{2}}$ with empty
intersection; they are also open in $K_{\alpha+1}$ and hence
$V_{x_{1}}$ and $V_{x_{2}}$ are open neighborhoods of $x_{1}$ and
$x_{2}$ respectively with empty intersection. \item[2)] There is
no $\vartheta\in\omega_{1}$ such that $x_{1},x_{2}\in
j_{\alpha+1}(\overline{N_{\vartheta}})$; therefore there are
$\vartheta_{1}, \vartheta_{2}\in\omega_{1}$ such that $x_{1}\in
j_{\alpha+1}(\overline{N_{\vartheta_{1}}})\simeq K_{\beta+1}$ with
$\beta+1<\alpha+1$ and $x_{2}\in
j_{\alpha+1}(\overline{N_{\vartheta_{2}}})\simeq K_{\gamma+1}$
with $\gamma+1<\alpha+1$. Now
$I=j_{\alpha+1}(\overline{N_{\vartheta_{1}}})\cap
j_{\alpha+1}(\overline{N_{\vartheta_{2}}})$ is a clopen subset of
$K_{\alpha+1}$ and hence
$j_{\alpha+1}(\overline{N_{\vartheta_{1}}})\backslash I$ and
$j_{\alpha+1}(\overline{N_{\vartheta_{2}}}) \backslash I$ are
disjoint open neighborhoods of $x_{1}$ and $x_{2}$ respectively.
\end{enumerate}
Therefore it turns out immediately that $K_{\alpha+1}$ is compact since
$j_{\alpha+1}$ is a continuous function from the compact space $\beta\omega$ to
the Hausdorff space $K_{\alpha+1}$.\finedim

Before proving the sequentiality of the space $K_{\alpha+1}$ we need to demonstrate that
properties S.4 and S.5 hold.

\begin{remark}\label{S4} In $K_{\alpha+1}$, if a nonconstant sequence $(x_{n})_{n\in\omega}$ of points
  $x_{n}\in{L_{\gamma_{n}+1}}$ with nondecreasing levels
   converges to a point $x\in{L_{\gamma+1}}$, then for the sequence
   $(\gamma_{n}+1)$ of ordinal numbers it holds that
   $\sup\{\gamma_{n}+1\}=\gamma$. (\emph{Properties S.4})
 \end{remark}
\iniziodim
  For every $\gamma+1<\alpha+1$ we apply the inductive hypothesis,
  since we have supposed that property S.4 holds in
   $K_{\gamma+1}$ for every $\gamma+1<\alpha+1$.\\
   Now we have to prove that for a non-constant sequence
of points $x_{n}\in{L_{\gamma_{n}+1}}$
    (where the sequence $(\gamma_{n}+1)$ is not decreasing)
    that converges to the point
   $x_{\infty}\in{L_{\alpha+1}}$ it holds that
   $\sup\{\gamma_{n}+1\}=\alpha$. Towards a contradiction, let us
   suppose that $\sup\{\gamma_{n}+1\}<\alpha$.
    In principle there are two different cases we have to analyse:
    \begin{itemize}
    \item[1)] from the sequence $(x_{n})_{n\in\omega}$ we can extract an
    injective subsequence $(x_{n_i})_{i\in\omega}$;
    \item[2)] from the sequence $(x_{n})_{n\in\omega}$ we can not extract any
    injective subsequence $(x_{n_i})_{i\in\omega}$.
    \end{itemize}
We can avoid considering the latter case: indeed, since $(x_{n})_{n\in\omega}$ is a
non-constant sequence, there are at least two points that appear infinite times and then
the sequence is not convergent to any point against the hypothesis.\\ In the former case
the sequence $(x_{n_i})_{i\in\omega}$ has to converge to $x_{\infty}$ too.
    If $\{x_{n_i}\}_{i\in\omega}$ was devoid of any accumulation point in
    $U_{\omega_{1}}/Q_{\omega_{1}}$,
     then by Remark \ref{ossdiscreto} it would be possible to find a countable infinity
     of pairwise disjoint clopen subsets of $\omega^{\ast}$; moreover these clopen subsets
      would satisfy the properties to be an $\omega_{1}$-family (notice that $\sup\{\gamma_{n}+1\}<\alpha$)
       and this would be inconsistent with Remark \ref{remark}.
       Then the subset $S=\{x_{n_i}:i\in\omega\}$
      has at least an accumulation point in $U_{\omega_{1}}/Q_{\omega_{1}}$;
       thus there exists a point $y\in U_{\omega_{1}}/Q_{\omega_{1}}$ (where certainly $l(y)=\delta+1<\alpha+1$)
        such that
        $y\in\overline{\{x_{n_i}:i\in\omega\}}$.
        Then let us consider an elementary neighborhood of $y$,
        $U_y$,
        which has to be homeomorphic to the space $K_{\delta+1}$;
        we can assert that infinite points of $S$ such that the supremum of their levels
         is equal to an ordinal number $\eta<\alpha$ are in $U_y$. We denote this set
         of points by $S'\subseteq S$; we know that property S.5 holds in
         $K_{\delta+1}$ and then from the injective sequence $S'$ it is possible to extract a sequence
         converging to a point of level $\eta+1<\alpha+1$.
          Therefore the sequence
       $(x_{n})_{n\in\omega}$ admits a subsequence which
        converges to a point of level strictly smaller than $\alpha+1$ against the hypothesis.\finedim

\begin{remark}\label{S5} In $K_{\alpha+1}$, from every injective sequence $S=(x_{n})_{n\in\omega}$
of points with nondecreasing levels such that
      $\sup\{l(x_{n})\}=\eta\leq\alpha$ it is possible to extract
       a subsequence converging to a point of level $\eta+1$.
       (\emph{Property S.5})
\end{remark}
\iniziodim If $\eta=0$ then the sequence $(x_{n})_{n\in\omega}$ is formed by
points of $\omega$; therefore there is an index $\gamma\in \omega_{1}$ such
that $|\{x_{n}\}_{n\in\omega}\cap N_{\gamma}|=\omega$: otherwise, if it turns
out that $|\{x_{n}\}_{n\in\omega}\cap N_{\gamma}|<\omega$ for every $\gamma\in
\omega_{1}$, then from $N_{\omega_{1}}=\{x_{n}\}_{n\in\omega}$ we are able to
construct an $\omega_{1}$-family and this is a contradiction. Then in
$\overline{N_{\gamma}}/\approx_{\alpha+1}$ there are infinite points of the
above sequence but $\overline{N_{\gamma}}/\approx_{\alpha+1}\simeq K_{\beta+1}$
with $\beta+1<\alpha+1$ and hence, since property S.5 holds in $K_{\beta+1}$ by
inductive hypothesis, it is possible to extract a subsequence
converging to a point of level $1$ from the starting sequence.\\
Suppose now that $0<\eta<\alpha$; let us choose an injective
subsequence $S'=(x_{n_i})_{i\in\omega}\subseteq S$ in such a way
that the sequence of the levels of the points converges upwards to
$\eta$; if $S'$ was devoid of any accumulation point in
$U_{\omega_{1}}/Q_{\omega_{1}}$, then by Remark \ref{ossdiscreto}
it would be possible to find a countable infinity
     of pairwise disjoint clopen subsets of $\omega^{\ast}$; moreover these clopen subsets
      would satisfy the properties to be an $\omega_{1}$-family (notice that $\sup\{l(x_{n_i})\}<\alpha$)
       and this would be inconsistent with Remark \ref{remark}. Thus $S'$ must have at least an accumulation point in
$U_{\omega_{1}}/Q_{\omega_{1}}$ and hence there exists a point $y\in
U_{\omega_{1}}/Q_{\omega_{1}}$ with $l(y)=\delta+1<\alpha+1$ such that
        $y\in\overline{\{x_{n_i}:i\in\omega\}}$. Then let us consider an elementary neighborhood of $y$,
        $U_y$,
        which has to be homeomorphic to the space $K_{\delta+1}$;
        we can assert that infinite points of the set $\{x_{n_i}:i\in\omega\}$ such that the limit and hence the supremum of their levels
         is equal to $\eta<\alpha$ are in $U_y$. We denote this set
         of points by $S''\subseteq S'$; we know that property S.5 holds in
         $K_{\delta+1}$ and then from the injective sequence $S''$ it is possible to extract a sequence
         converging to a point of level $\eta+1<\alpha+1$.\\
 If $\eta=\alpha$ then let us choose again an injective
subsequence $S'=(x_{n_i})_{i\in\omega}\subseteq S$ in such a way that the levels of the
points $x_{n_i}$ converges upwards to $\alpha$; the sequence $S'$ converges to
    $x_{\infty}$, since its points fall eventually in every neighborhood of $x_{\infty}$.\finedim

Now we are able to prove the sequentiality of $K_{\alpha+1}$.
 \begin{lemma} $K_{\alpha+1}$ is sequential.
 \end{lemma}
\iniziodim Let us begin by proving that $B_{\alpha+1}=K_{\alpha+1}\backslash
\{x_{\infty}\}$ is sequential, i.e. by showing that if $F$ is a sequentially
closed subset of $B_{\alpha+1}$ then it is closed. Let us suppose that $F$ is
sequentially closed and let us show that for every $x\in B_{\alpha+1}\backslash
F$ there exists an elementary neighborhood $\hat{U}_{x}$ of $x$ such that
$\hat{U}_{x}\subseteq B_{\alpha+1}\backslash F$. If $x\in
B_{\alpha+1}\backslash F$, then there exists an open neighborhood of $x$,
$U_{x}\subseteq B_{\alpha+1}$ with the peculiarity that $U_{x}\simeq
K_{\beta+1}$ (with $\beta+1<\alpha+1$) which is a compact sequential space.
Notice that $x\notin F\cap U_{x}$; if $F\cap U_{x}=\emptyset$, then $U_{x}$ is
an elementary neighborhood containing $x$ such that $U_{x}\subseteq
B_{\alpha+1}\backslash F$. If instead $F\cap U_{x}\neq\emptyset$, since $F$ is
sequentially closed in $B_{\alpha+1}$, then $F\cap U_{x}$ is sequentially
closed in $U_{x}$ (otherwise, if $F\cap U_{x}$ is not sequentially closed in
$U_{x}$, hence we have a sequence in $F\cap U_{x}$ with its limit point in
$U_{x}\backslash F$; we can see this sequence as a sequence in $F$ with its
limit point out of $F$ and then $F$ is not sequentially closed against the
hypothesis). It follows that $F\cap U_{x}$ is closed in $U_{x}$ since $U_{x}$
is sequential and hence it is compact; let us consider the open cover of $F\cap
U_{x}$ formed by elementary neighborhoods of points in $F\cap U_{x}$ not
containing $x$. From this open cover it is possible to extract a finite
subcover $\bigcup_{i=1}^{\overline{n}}U_{y_{i}}\supseteq F\cap U_{x}$. Then
$\hat{U}_{x}= U_{x}\backslash \bigcup_{i=1}^{\overline{n}}U_{y_{i}}$ is an open
neighborhood of $x$ which has empty intersection with $F$. Thus we can conclude
that $B_{\alpha+1}$ is sequential.\\ Now we have still to demonstrate that, if
$F$ is sequentially closed in $K_{\alpha+1}$ and $x_{\infty}\notin F$, then
$x_{\infty}\notin\overline{F}$. Towards a contradiction, suppose that
$x_{\infty}\notin F$ and, at the same time, $x_{\infty}\in\overline{F}$. Since
$x_{\infty}\notin F$ then either $F$ is finite (and in this case the point
$x_{\infty}\notin\overline{F}$ against the hypothesis) or $F$ is infinite and
in this second case from $F$ it is not possible to extract any injective
sequence of points with nondecreasing levels such that the supremum of the
levels is equal to $\alpha$; indeed if such a sequence existed, by Remark
\ref{S5} from this sequence it would possible to extract a subsequence
converging to $x_{\infty}$ and then $x_{\infty}$ would stay in $F$ (since $F$
is sequentially closed) against the hypothesis. Now if $\alpha$ is a successor
ordinal, there exists at most a finite number of points of level
$\alpha=\gamma_{0}$ in $F$ that we call $z_{1},z_{2},\ldots,z_{m}$; let us
consider an elementary neighborhood $U_{z_{i}}$ for each of these points and
let us set $G_{1}=F\backslash \bigcup_{i=1}^{m}U_{z_{i}}\subseteq F$. We assert
that either $G_{1}$ is finite (and in this case it turns out that
$x_{\infty}\notin\overline{F}$ against the hypothesis) or $G_{1}$ is infinite
and in this second case from $G_{1}$ it is not possible to extract any
injective sequence of points with nondecreasing levels such that the supremum
of the levels is equal to $\alpha-1=\gamma_{1}$; indeed if such a sequence
existed, by Remark \ref{S5} from this sequence it would possible to extract a
subsequence converging to a point of level $\alpha$ different from
$z_{1},z_{2},\ldots,z_{m}$ and then also this point would stay in $F$ against
our assumption. If instead $\alpha$ is a limit ordinal, it is not true that for
every $\gamma\in\alpha$ there exists $x\in F$ such that $l(x)>\gamma+1$
(otherwise $x_{\infty}\in F$ which is sequentially closed) and hence there
exists an index $\overline{\gamma}\in\alpha$ such that for every $x\in F$ it
turns out that $l(x)\leq\overline{\gamma}+1<\alpha$. Therefore we can assert
that in $F$ there are at most a finite number of elements of level
$\overline{\gamma}+1$ that we call $y_{1},y_{2},\ldots,y_{k}$: indeed if we had
an infinite number of these points, it would possible to extract a subsequence
converging to a point of level $(\overline{\gamma}+1)+1$ and this point would
stay again in $F$ but this is against what we have just remarked. Let us
consider an elementary neighborhood $U_{y_{i}}$ for each of these points and
let us call $G_{1}=F\backslash \bigcup_{i=1}^{k}U_{y_{i}}\subseteq F$. We can
say that either $G_{1}$ is finite (and in this case the point
$x_{\infty}\notin\overline{F}$) or $G_{1}$ is infinite and in this second case
from $G_{1}$ it is not possible to extract any injective sequence of points
with nondecreasing levels such that the supremum of the levels is equal to
$\gamma_{1}=\overline{\gamma}<\alpha$; indeed if such a sequence existed, by
Remark \ref{S5} from this sequence it would be possible to extract a
subsequence converging to a point of level $\overline{\gamma}+1$ different from
$y_{1},y_{2},\ldots,y_{k}$ and then also this point would stay in $F$ against
what we have
assumed.\\
 In each case we have constructed a sequentially closed subset
  $G_{1}$ from which it is not possible to extract any injective
sequence of points with nondecreasing levels such that the supremum of the
levels is equal to $\gamma_{1}<\alpha$; moreover $G_{1}$ is the complement of a
finite number of elementary neighborhoods in $F$. Then it is possible to repeat
the procedure and to find step by step a decreasing sequence of ordinals
$\gamma_{0}>\gamma_{1}>\gamma_{2}>\ldots>\gamma_{n}>\ldots$ and corresponding
subsets $G_{1}\supseteq G_{2}\supseteq\ldots\supseteq G_{n}\supseteq \ldots$.
This sequence has to be finite and then we find a finite set $G_{\overline{n}}$
after a finite number of steps. Trivially we can cover $G_{\overline{n}}$ by a
finite number of elementary neighborhoods; moreover $G_{\overline{n}}$ has been
constructed as complement of a finite number of elementary neighborhoods in
$F$. Then it turns out that it is possible to cover $F$ with finitely many
elementary neighborhoods of points of level smaller than $\alpha+1$ and hence
it follows that $x_{\infty}\notin \overline{F}$. A contradiction. \finedim

Now we want to show that every point in $K_{\alpha+1}$ belongs to
the closure of $L_{0}$, i.e that the set $L_{0}$ is dense in
$K_{\alpha+1}$.
\begin{remark}\label{adalphasuc} For every $x\in
K_{\alpha+1}$ it holds that $x\in \overline{L_{0}}$, i.e.
$\overline{L_{0}}=K_{\alpha+1}$.
\end{remark}
\iniziodim Let $x$ be a point in $K_{\alpha+1}$ with $l(x)=\beta+1<\alpha+1$
and let $V$ be a non-empty neighborhood of $x$; then there exists an open
elementary neighborhood $U_{x}\simeq K_{\beta+1}\subseteq V$ and it turns out
that $j_{\alpha+1}^{-1}(U_{x})$ is a non-empty open subset in $\beta\omega$.
Therefore there exists a free or a fixed ultrafilter $\mathcal{U}$ such that
$\mathcal{U}\in j_{\alpha+1}^{-1}(U_{x})$. If $\mathcal{U}$ is fixed we
trivially finish; if $\mathcal{U}$ is a free ultrafilter, since
$j_{\alpha+1}^{-1}(U_{x})$ is an open subset, there is an infinite subset $U'$
of $\omega$ with $U'\in{\mathcal{U}}$ such that $(U')^{\ast}\cup U'\subseteq
j_{\alpha+1}^{-1}(U_{x})$; then $U'\subseteq\omega$
  (with $|U'|=\omega$) is such that $U'\subseteq j_{\alpha+1}^{-1}(U_{x})$
  and hence $W=j_{\alpha+1}(U')\subseteq U_{x}$; we can conclude that $U_{x}\cap
   L_{0}\supseteq W\cap L_{0} \neq\emptyset$.\\
Now let us consider $x_{\infty}\in K_{\alpha+1}$ and let $U$ be an open
neighborhood of $x_{\infty}$. By Lemma \ref{lemma intorni punti con alpha
successore} there exists an open subset $A_{x_{\infty}}=K_{\alpha+1}\backslash
\bigcup U_{x}\subseteq U$; therefore $j_{\alpha+1}^{-1}(A_{x_{\infty}})$ is a
non-empty open subset of $\beta\omega$ and hence we can proceed as above.
\finedim

Since we have proved that the space $K_{\alpha+1}$ is sequential
and that $\overline{L_{0}}=K_{\alpha+1}$, Remark \ref{S4} allows
us to conclude that the level of each point is larger or equal to
its order of sequentiality with respect to $L_{0}$. We have to
prove a last remark before concluding that the level of each point
is exactly equal to its sequential order with respect to the set
$L_{0}$.
\begin{remark}\label{5)alphasuc}
Let $A$ be a closed subset in $\bigcup_{\gamma+1\leq\eta}L_{\gamma+1}$ with
$\eta\leq\alpha$. Then it follows that $\overline{A}\cap
\bigcup_{\gamma+1\leq\eta+1}L_{\gamma+1}=seqcl(A)$.
\end{remark}
\iniziodim Since $A$ is closed in $\bigcup_{\gamma+1\leq\eta}L_{\gamma+1}$,
then $A$ is sequentially closed in $\bigcup_{\gamma+1\leq\eta}L_{\gamma+1}$ and
hence there is no sequence in $A$ converging to some point of
$\bigcup_{\gamma+1\leq\eta}L_{\gamma+1}\backslash A$. Notice that by Remark
\ref{S4} it turns out that $seqcl(A)\subseteq\overline{A}\cap
\bigcup_{\gamma+1\leq\eta+1}L_{\gamma+1}$. We want to prove that
$\overline{A}\cap \bigcup_{\gamma+1\leq\eta+1}L_{\gamma+1}\subseteq seqcl(A)$.
Let $x$ be a point in $(\overline{A}\backslash
A)\cap(\bigcup_{\gamma+1\leq\eta+1}L_{\gamma+1})$; trivially it holds that
$l(x)=\eta+1$. Since $K_{\alpha+1}$ is sequential and $x\in\overline{A}$, it
turns out that there exists an index $\beta\in\omega_{1}$ such that $x\in
seqcl_{\beta}(A)$; we state that $\beta=1$. Towards a contradiction, let us
suppose that $x\notin seqcl_{1}(A)$, i.e. let us suppose that no sequence in
$A$ converges to $x$. Then $x$ is the limit point of a sequence whose elements
are in some sequential closure of $A$ and not in $A$, i.e. $x$ is the limit
point of a sequence $(y_{\beta_{i}+1})_{i\in\omega}$ with
$\sup\{\beta_{i}+1\}\geq\eta+1$ but this is absurd since $l(x)=\eta+1$: indeed
if it was correct, in $K_{\alpha+1}$ there would exist a sequence
$(y_{\beta_{i}+1})_{i\in\omega}$ with $\sup\{\beta_i+1\}\neq\eta$ converging to
a point of level $\eta+1$ and this is inconsistent with Remark \ref{S4}.
\finedim

Finally we can prove the following crucial lemma.
\begin{lemma}\label{lemordalphasuc} In $K_{\alpha+1}$ the order of sequentiality of a point
of level $\beta+1$ with respect to $L_{0}$ is $\beta+1$ and
$K_{\alpha+1}$ is a space with sequential order $\alpha+1$.
\end{lemma}
\iniziodim Notice that the points of level $0$ and $1$ have sequential order
respectively $0$ and $1$ with respect to the set $L_{0}$. Now consider the set
$\bigcup_{\gamma+1\leq\eta}L_{\gamma+1}$ with $\eta\leq\alpha$; it complies
with the hypotheses of Remark \ref{5)alphasuc} since it is closed in
$\bigcup_{\gamma+1\leq\eta}L_{\gamma+1}$ and hence it holds that
$$\overline{\bigcup_{\gamma+1\leq\eta}L_{\gamma+1}}\bigcap
\bigcup_{\gamma+1\leq\eta+1}L_{\gamma+1}=$$$$\{y \in
\bigcup_{\gamma+1\leq\eta+1}L_{\gamma+1}: \forall U_{y}, (U_{y}\cap
\bigcup_{\gamma+1\leq\eta}L_{\gamma+1})\neq\emptyset\}=$$$$
\bigcup_{\gamma+1\leq\eta+1}L_{\gamma+1}=seqcl(\bigcup_{\gamma+1\leq\eta}L_{\gamma+1}).$$
This result together with Remark \ref{adalphasuc} allows us to conclude that
the level of each point is smaller or equal to its order of sequentiality with
respect to the set $L_{0}$. But we have already remarked that the level of each
point is larger or equal to its order of sequentiality with respect to the set
$L_{0}$ and hence we can conclude that the level of each point is exactly equal
to its order of sequentiality with respect to the set $L_{0}$.\\ Then the space
$K_{\alpha+1}$ has sequential order equal to $\alpha+1$, since $x_{\infty}\in
\overline{L_{0}}$ and $x_{\infty}$ has sequential order equal to $\alpha+1$
with respect to $L_{0}$. \finedim

Let us finally check that properties $S.1$ to $S.6$ hold in the
space $K_{\alpha+1}$.
\begin{enumerate}
\item [S.1] The space $K_{\alpha+1}$ can be uniquely represented
in the form of $$K_{\alpha+1}=L_{0}\bigsqcup (\bigsqcup_{\gamma\leq\alpha}
L_{\gamma+1}).$$
 The points of level
$\gamma+1$ with $\gamma\in[0,\alpha]$, i.e. the points belonging
     to the set $L_{\gamma+1}$, have sequential order equal to
    $\gamma+1$ with respect to $L_{0}$: see Lemma \ref{lemordalphasuc}.
 \item [S.2] The set
$L_{\alpha+1}$ consists of the unique point $x_{\infty}$.
 \item [S.3] Every point in $K_{\alpha+1}$ of nonzero level has a basis formed by clopen
     subsets called elementary; if $U$ is an elementary neighborhood of a point of level
     $\gamma+1$, then the relation $\approx_{\alpha+1}$ restricted to
    $\widetilde{U}=j^{-1}_{\alpha+1}(U)$ produces a compact space homeomorphic to $K_{\gamma+1}$: for the points of level
     smaller than $\alpha+1$, S.3 is true by inductive hypothesis
      while for $x_{\infty}$ the property is correct since each of its elementary neighborhood is homeomorphic to the whole
        space $K_{\alpha+1}$ (see Lemma
      \ref{lemma intorni punti con alpha successore}).
  \item [S.4] For every $\gamma\leq\alpha$, if a nonconstant sequence $(x_n)_{n\in\omega}$ of points
  $x_{n}\in{L_{\gamma_{n}+1}}$, with nondecreasing levels,
   converges to a point $x\in{L_{\gamma+1}}$, then for the sequence $(\gamma_{n}+1)_{n\in\omega}$
    of ordinal numbers it holds that
    $\sup\{\gamma_{n}+1\}=\gamma$: see Remark \ref{S4}.
    \item [S.5] For every $\gamma\leq\alpha$, from every injective sequence $(x_{n})_{n\in\omega}$ of points
        $x_{n}\in{L_{\gamma_{n}+1}}$ with nondecreasing levels such that
      $\sup_{n\in\omega}\{\gamma_{n}+1\}=\gamma$, it is possible to extract
       a subsequence converging to a point of level $\gamma+1$: see Remark \ref{S5}.
  \item [S.6] If $\left\{N_{i}\right\}_{i\in\omega}$ is a countable family
       of pairwise disjoint infinite subsets $N_{i}$
     of $\omega$ and if it holds that for every $i\in\omega$ a relation of type $\beta_{i}+1$ is given on
     $\overline{N_{i}}$ in such a way that the sequence of ordinals $(\beta_{i}+1)_{i\in\omega}$ is not
     decreasing and $\sup\left\{\beta_{i}+1\right\}=\alpha$, then it is possible
     to extend the relation obtained on $\bigcup_{i=1}^{\infty}\overline{N_{i}}$ to a relation of $\beta\omega$ of type
     $\alpha+1$: we have just constructed it.
\end{enumerate}
\begin{remark}
Notice that every Ba\v{s}kirov's space of sequential order a successor ordinal is a
scattered space such that the sequential order of each point is equal to its scattering
level.
\end{remark}
 Finally we
can state the following theorem.
\begin{theorem}(CH) Let $\alpha$ be any ordinal less than or equal
to $\omega_1$. There exists a compact sequential Hausdorff quotient space of
$\beta\omega$ with sequential order $\alpha$.
\end{theorem}
\hspace{-0.5cm} \textbf{Acknowledgments.} I wish to thank Gino Tironi,
Professor at the University of Trieste, for giving me the opportunity to work
on this challenging topic, raising my
interest in topology and for patiently answering all my questions.\\
I am grateful to Camillo Costantini, Assistant Professor at the University of
Torino, for helping me to get more familiarity with some techniques which are
of use to investigate the notions related to my research field.

\bibliographystyle{plain}

\begin{thebibliography}{10}

\smallskip

\bibitem{arhangelskii} A. V. Arhangel'skii and S. P. Franklin, \textit{Ordinal invariants for topological spaces},
Michigan Math. J. \textbf{15} (1968), 313--320.

\smallskip

\bibitem{balogh} Zolt\'{a}n T. Balogh, \textit{On compact Hausdorff spaces of countable tightness},
Proc. Amer. MAth. Soc., \textbf{105(3)} (1989), 755--764.

\smallskip

\bibitem{baskirov} A. I. Ba{\v{s}}kirov, \textit{The classification of quotient maps and sequential
bicompacta}, Soviet Math. Dokl. \textbf{15} (1974), 1104--1109.

\smallskip

\bibitem{dow1} Alan Dow, \textit{On MAD families and sequential
order},\\
http://www.math.uncc.edu/${}_{\textrm{\symbol{126}}}$adow/Others.htmal
under the title ``On the sequential order of Compact spaces''.

\smallskip

\bibitem{dow2} Alan Dow, \textit{Sequential order under MA},
 Topology Appl. \textbf{146-147} (2005), 501--510.

\smallskip

\bibitem{engelking} Ryszard Engelking, \textit{General Topology},
Translated from the Polish by the author. 2nd ed. Sigma Series in
Pure Mathematics, 6. Berlin: Helderman Verlag, 1989.

\bibitem{franklin1} S.P. Franklin, \textit{Spaces in which sequences suffice},
Fund. Math. \textbf{57} (1965), 107--115.

\smallskip

\bibitem{franklin2} S.P. Franklin, \textit{Spaces in which sequences suffice. II},
Fund. Math. \textbf{61} (1967), 51--56.

 \smallskip

\bibitem{kannan} V. Kannan, \textit{Ordinal invariants in topology II. Sequential order of compactifications},
 Compositio Math. \textbf{39(2)} (1979), 247--262.

\smallskip

\bibitem{kannan2} V. Kannan, \textit{Ordinal invariants in topology III. Simplest compact spaces for sequential order},
 Math. Today \textbf{1} (1983), 65-80.

\smallskip

\bibitem{vanmill} Jan van Mill, \textit{An introduction to $\beta\omega$},
In Handbook of set-theoretic topology, North- Holland, Amsterdam,
(1984), 503-567.

 \smallskip

\bibitem{rudin} Walter Rudin, \textit{Homogeneity problems in the thoery of \v{C}ech compactifications},
Duke Math. J \textbf{23} (1956), 409-419.

\end{thebibliography}

\end{document}